\input amstex
\documentstyle{amsppt}

\input epsf.tex
\magnification=1200

\pagewidth{6.2truein}
\pageheight{9truein}

\baselineskip=18pt
\NoRunningHeads

\def\figure#1#2#3{\vskip6pt\epsfysize#1truein
                \centerline{\epsfbox{#2}}
                \medskip
                \centerline{ {\smc Figure} #3}\vskip6pt}

\topmatter
\title Meridional surfaces and $(1,1)$-knots
\endtitle
\author Mario Eudave-Mu\~noz and Enrique Ram\'{\i}rez-Losada
\endauthor

\address Instituto de Matem\'aticas, UNAM,\\ Ciudad Universitaria,
04510 M\'exico D.F., MEXICO \\
E-mail: mario\@matem.unam.mx\\kikis\@cimat.mx \endaddress

\abstract
We determine all $(1,1)$-knots which admit an essential meridional surface,
namely, we give a construction which produces $(1,1)$-knots having essential
meridional surfaces, and show that if a $(1,1)$-knot admits an essential
meridional surface then it comes from the given construction.

\endabstract


\endtopmatter

\document

\medskip
\head {1. Introduction} \endhead
\medskip

We are interested in finding closed incompressible surfaces or meridional
incompressible surfaces in the complement
of knots in $S^3$. This problem has been studied for many classes of knots,
for example for 2-bridge knots [GL],[HT], Montesinos knots [O], alternating knots
[M], closed 3-braids [LP],[F], etc.

Let $k$ be a knot in a 3-manifold $M$.
Recall that a surface properly embedded in the exterior of $k$ in $M$
is meridional if its boundary consists of a collection of meridians of $k$
(and possibly some curves on $\partial M$).
A surface $S$ properly embedded in the exterior of $k$ (either closed or meridional),
is meridionally compressible if there is a disk $D$ embedded in $M$, so that
$D\cap S = \partial D$ is an essential curve in $S$ (non-contractible in $S$ and
non-parallel in $S$ to a meridian of $k$), such that $k$ intersects $D$
transversely in exactly one point (or in other words, there is an annulus $A$
embedded in the exterior of $k$, with $\partial A= \partial_0 A \cup \partial_1 A$,
such that $A\cap S = \partial_0 A$ and $\partial_1 A$ is a meridian of $k$); if there
is no such a disk, we say that the surface is meridionally incompressible. We say
that a meridional surface $S$ is essential if it is incompressible and
meridionally incompressible.

A class of knots that has been widely studied is that of tunnel number one
knots (a knot $k$ has tunnel number one if there is an arc $\tau$ embedded in
$S^3$, with $k\cap \tau =\partial \tau$ such that $S^3-int\, \eta (k\cup \tau)$ is a
genus 2 handlebody). By work of Gordon and Reid [GR], these knots do not admit
any planar meridional incompressible surface. However, some of these knots do
admit closed incompressible surfaces or meridional incompressible surfaces of
genus one or greater. Morimoto and Sakuma [MS] determined all tunnel number one
knots whose complements contain an incompressible torus. The first author
constructed for any $g\geq 2$ infinitely many tunnel number one knots whose
complements contain closed incompressible surfaces of genus $g$ [E1][E2], and for any
$g\geq 1$ and $h\geq 1$, infinitely many tunnel number one knots
whose complements contain an essential  meridional surface of genus $g$ and $2h$
boundary components [E2]. An interesting problem would be to determine,
characterize  or classify all tunnel number one knots which admit a closed or
meridional incompressible surface. With this objective we initiated studying
incompressible surfaces in a special class of tunnel number one knots, that of
$(1,1)$-knots. A knot $k$ is said to be a $(1,1)$-knot if there is a standard
torus $T$ in $S^3$ such that $k$ is a 1-bridge knot with respect to $T$, that is,
$T$ and $k$ intersect transversely in two points, which divide $k$ into two arcs,
and such that each arc is isotopic to an arc lying on $T$.  Equivalently,
consider $T\times I\subset S^3$,  where $I=[0,1]$; $k$ is a $(1,1)$-knot if
$k$ lies in $T\times I$, so that $k\cap T\times \{ 0\} = k_0$ is an arc,
$k\cap T\times \{ 1\} = k_1$ is an arc, and $k\cap T\times (0,1)$ consists
of two straight arcs, i.e., arcs which cross every torus
$T\times \{x\}$ transversely. The family of $(1,1)$-knots is an
interesting class of knots which contains torus knots, 2-bridge knots, and which
has received much attention recently (see for example [CM],[CK], [GMM]).
All satellite tunnel number one knots classified by
Morimoto and Sakuma are in fact $(1,1)$-knots, and many of the knots constructed
in [E1] are also $(1,1)$-knots. In the paper [E3] it is proved that if $k$ is a
$(1,1)$-knot whose complement contains a closed incompressible and meridionally
incompressible surface, then $k$ comes from the construction given in [E1].

The class of $(1,1)$-knots can be defined for any manifold $M$ admitting a Heegaard
decomposition of genus 1, i.e., it can be defined when $M= S^3,\ S^1\times S^2$, or a
lens space $L(p,q)$, the definition is the same, with $T$ being a Heegaard torus.

Let $M= S^3,\ S^1\times S^2$, or a lens space $L(p,q)$.
In this paper we construct all $(1,1)$-knots in $M$ having a
meridional essential surface. First, we give a general construction which
produces $(1,1)$-knots which admit a meridional essential surface (Section 2),
and then prove that if a $(1,1)$-knot admit a meridional essential surface then
it comes from the given construction (Section 3). In particular we show that
for any given integers $g\geq 1$ and $h\geq 0$, there exist $(1,1)$-knots which
admit a meridional essential surface of genus $g$ and $2h$ boundary components.
As we said before, $(1,1)$-knots in $S^3$ do not admit any meridional essential
surface of genus $0$, but when
$M=L(p,q)$, we show that for any given integer $h\geq 1$, there exist $(1,1)$-knots
which admit a meridional essential surface of genus $0$ and $2h$ boundary
components.

We remark that the knots constructed here are
different from the ones in [E2]; in fact, most of the knots
constructed in [E2] do not seem to be $(1,1)$-knots.

Throughout, 3-manifolds and surfaces are assumed to be compact and orientable. If $X$
is contained in a 3-manifold $M$, then $\eta (X)$ denotes a regular neighborhood of
$X$ in $M$.

\medskip
\medskip
\head {2. Construction of the meridional surface} \endhead

In this section we construct meridional surfaces for $(1,1)$-knots. We construct
the surfaces by pieces, that is, we define some surfaces in a product $T\times I$
or in a solid torus, getting 6 types of basic pieces, denoted by
$\Cal {\tilde A}$, $\Cal {\tilde B}$, $\Cal {\tilde C}$, $\Cal {\tilde D}$,
$\Cal {\tilde E}$, $\Cal {\tilde F}$.
By assembling the pieces appropriately we get a $(1,1)$-knot
having an essential meridional surface.

\medskip

\noindent {\bf 2.1} Let $T$ be a torus, and for any two real numbers $a,b$, $a < b$,
define
$N_{a,b}=T\times [a,b]$, and $\partial_a N_{a,b}=T\times \{a\}$, $\partial_b
N_{a,b}=T\times \{b\}$. A properly embedded arc in $N_{a,b}$ is straight if it
intersects each torus $T\times \{x\}$ in one point, for all $a\leq x \leq b$.
Let $A$ be a once punctured annulus (i.e. a pair of pants) properly
embedded in $N_{a,b}$, such that a boundary component of $A$, say
$\partial_0 A$, and which we call the puncture, lies on $\partial_a N_{a,b}$
and is a trivial curve in this torus; the other components, denoted by
$\partial_1 A$, are a pair of essential curves on  $\partial_b N_{a,b}$.
Assume that $A$ has been isotoped so that $A$ has only one saddle singularity with
respect to the projection $N_{a,b} \rightarrow [a,b]$, that is, there is a  real
number $y$,
$a< y < b$, so that $A\cap  (T\times \{ z\})$ consists of a trivial curve in
$T\times \{ z\}$ for $a\leq z < y$, $A\cap  (T\times \{ z\})$ consists of two essential
curves in  $T\times \{ z\}$ for $y < z\leq  b$, and $A\cap  (T\times \{ z\})$ consists
of a curve with a selfintersection if $z=y$. Two such punctured annuli $A_1$, $A_2$
are parallel if $A_1\times I$ is embedded in $N_{a,b}$ such that
$A_1\times \{ 0 \}=A_1$, $A_1\times \{ 1\}=A_2$, and
$\partial_0 A_1\times I $  ($\partial_1 A_1\times I $) lies on $\partial_a N_{a,b}$
($\partial_b N_{a,b}$).

We defined in the introduction the notion of meridional incompressibility of a surface
with respect to a knot in a 3-manifold, but note that the same definition can be given
with respect to a link or with respect to a collection of arcs properly embedded in a
3-manifold.

\medskip
\subhead {2.2 Pieces of type $ \Cal {\tilde A}$} \endsubhead
\medskip

Let $N_{a,c}=N_{a,b}\cup N_{b,c}$, for $a < b < c$.
Denote by $A_1,\dots, A_r$ a collection of
disjoint, parallel, once punctured annuli properly embedded in $N_{a,b}$ as in 2.1.
Note that the curves $\partial_0 A_i$ are a collection of $r$  trivial,
nested curves in $\partial_a N$, i.e., each of them bounds a disk $D_i$,
such that $D_i\subset D_{i+1}$, for
$i\in \{ 1,2,...,r-1\}$. The curves $\partial_1A_i$ are a collection of $2r$
parallel essential curves in  $\partial_b N$. Let $A_1',...,A_r'$ be another
collection of properly embedded parallel once punctured annuli in $N_{b,c}$,
again as in 2.1, but such that the curves
$\partial_0 A_i'$, i.e. the punctures, lie on  $\partial_c N_{b,c}$
and bound disks $D_i'$, and the
curves $\partial_1 A_i'$ lie on $\partial_b N_{b,c}$. Suppose furthermore that
$\{ \partial_1 A_i \} = \{ \partial_1 A_i' \}$. Note that the boundaries of an
annulus in one side may be identified to curves belonging to two different annuli
in the other side.
Let $\Cal A'$ be the union of the annuli $A_i$ and $A_i'$. $\Cal A'$ is a surface
properly embedded in $N_{a,c}$. Note that each component of $\Cal A'$ is a torus
with an even number of punctures, and that $\Cal A'$ could be connected. Let
$\eta(\partial_b N_{a,b})$ be a regular neighborhood of $\partial_b N_{a,b}$
in $N_{a,c}$ such that $\eta (\partial_b N_{a,b})\cap {\Cal A'}$ consists of a
collection of $2r$ parallel vertical annuli. Let $t_{a,c}$ be two straight
arcs properly embedded in $N_{a,c}$ such that
$t_{a,c}\cap \partial_a N_{a,c}$ consists of two points contained in $D_1$,
and similarly $t_{a,c}\cap \partial_c N_{a,c}$ consists of two points
contained in $D_1'$. Suppose also that the arcs of $t_{a,c}$ intersect
$\Cal A'$ in finitely
many points, all of which lie in $\eta(\partial_b N_{a,b})$. Assume these
intersection points lie at different heights; denote the points and the height at
which they occur by $x_1,x_2,\dots,x_n$, which are ordered according to their heights.
It should be clear from the context whether we refer to an intersection point or to a
level in which the point lies.

\medskip
\noindent {\bf 2.2.1} Let $\Cal A = \Cal A'- int\, \eta (t_{a,c})$.
We call a product $N_{a,c}$ together with a surface $\Cal A$ and a
pair of arcs $t_{a,c}$, a piece of type $\tilde \Cal A$. We assume that
a piece of type $\tilde \Cal A$ satisfies the following:
\roster

\item The part of $A_1$ up to level $x_1$, $D_1$ and an annulus $E_1$ in
$T\times \{x_1 \}$ bound a solid torus $N_1\subset N_{a,x_1}$. Suppose there is
no a meridian disk $D$ of $N_1$, disjoint from $t_{a,c}$, and whose boundary
consists of an arc on $A_1$, and one arc in $E_1$.

\noindent Note that if this condition is not satisfied, and so there
is such a meridian disk, then the subarc $t$ of $t_{a,c}$ which
starts at $D_1$ and ends at $x_1$, is parallel onto $A_1$, that is,
there is a disk $F$ embedded in $N_1$, with $\partial F = t\cup
\alpha\cup \beta$, where $\alpha$ is an arc on $D_1$ and $\beta$ an
arc on $A_1$, and $int F\cap t_{a,c}=\emptyset$. So we can slide the
arc $t$ until is at the level of the disk $D_1$, and so there is a
meridional compression disk for $\Cal A$. That is, if the condition
is not satisfied, then the surface $\Cal A$ will be meridionally
compressible. On the other hand, it is not difficult to see that if
the arc $t$ is parallel to an arc on $A_1$, then there is a meridian
disk for $N_1$ disjoint from $t_{a,c}$.

\item One arc of $t_{a,c}$ start at $D_1$ and arrives to the point $x_1$. Let
$t'$ be the other arc of $t_{a,c}$. Assume that either (i) $t'$ is disjoint
from $\Cal A'$; or (ii) $x_1$ and the first point of intersection of
$t'$ with $\Cal A'$, say $x_j$, lie on different vertical annuli
of $\eta(\partial_b N_{a,b})\cap\Cal A'$; or (iii) $x_j$ and $x_1$ lie on the
same vertical annulus, but there is no disk
$D \subset \eta(\partial_b N_{a,b})$, so that
$\partial D =\alpha \cup \beta\cup\gamma$, where  $\alpha\subset t'$,
$\beta \subset \Cal A'$, which connects $x_j$ with a point $x_1'$ at level
$T\times \{x_1\}$, $\gamma$ is an arc on
$T\times \{x_1\}$, and $int\, D$ is disjoint from both $t_{a,c}$ and $\Cal A'$. In
other words, $t'$ cannot be slided to lie at level $T\times \{x_1\}$.

\noindent In particular, if the first point of intersection of $t'$ with $\Cal A'$ is $x_2$,
then we are assuming that $x_1$ and $x_2$ lie on different vertical annuli.
Note that if this condition is not satisfied then condition (1) would be meaningless,
i.e., we could always slide the arcs and find a meridional compression disk for $\Cal A$.

\item Suppose a subarc of $t_{a,c}$ intersects a vertical annulus in two points,
say $x_i$, $x_j$, but does not intersect any other annulus between these two points.
Then this subarc, say $\beta$, is not parallel to an arc
on the vertical annulus. That is, if $D$ is any disk in $N_{a,c}$ with
interior disjoint from $\Cal A'$, such that $\partial D=\alpha \cup \beta$,
where $\alpha$ is an arc on the vertical annulus,
then $t_{a,c}$ necessarily intersects $int\, D$. If this condition is not
satisfied then the surface $\Cal A$ will be clearly compressible.

\item For the part of $A_1'$ between levels $T\times \{x_n\}$ and
$T\times \{c\}$ we have a similar condition as in (1).

\item For $D_1'$ and the arcs adjacent to it we have a similar condition as in (2).

\item  In case $t_{a,c}$ is disjoint from $\Cal A'$, we assume the
following:
In this case the annulus $A_i$ is glued to the annulus $A_i'$. Let $k$ be
the knot obtained by joining the endpoints of $t_{a,c}$ contained in $D_1$
with an arc lying on $D_1$, and by joining the endpoints contained on $D_1'$
with an arc lying on $D_1'$. Then we get a knot contained in the solid torus
bounded by $D_1 \cup A_1 \cup A_1'\cup D_1'$.
Assume this knot has wrapping number $\geq 2 $ in such
a solid torus. This is required to avoid meridional compression disks.

\endroster

Note that it is not difficult to construct
pieces of type $\tilde \Cal A$ where all of these conditions are satisfied.
See Figure 1.

\remark {Remark} Condition 2.2.1(1) has an alternative description.
Embed the solid torus $N_1$
in $S^3$ so that it is an standard solid torus and $A_1\cap (T\times\{x_1 \})$
is a preferred longitude of such solid torus. Connect the endpoints of
$t_{a,c}\cap N_1$ contained in $D_1$ with an arc lying in $D_1$. The other
endpoints of $t_{a,c}\cap N_1$ lie on $E_1$; let $\alpha$ be the boundary of a
meridian disk of $N_1$, passing through these points and which is disjoint
from $D_1$. The endpoints of $t_{a,c}\cap N_1$ separate $\alpha$ into two arcs;
joint these points by the arc of $\alpha$ that is not contained in $E_1$.
This defines a knot $k$ in
$S^3$, which by construction has a presentation with two maxima.
Note that if condition 2.2.1(1) is not satisfied, then $k$ will be the
trivial knot. So, a sufficient condition for the condition 2.2.1(1) to be satisfied
is that $k$ is a non-trivial 2-bridge knot. With a little work, it can be shown that
this condition is also necessary. Observe that any 2-bridge knot can result
from this construction, as shown in Figure 3.


\midinsert
\figure{3}{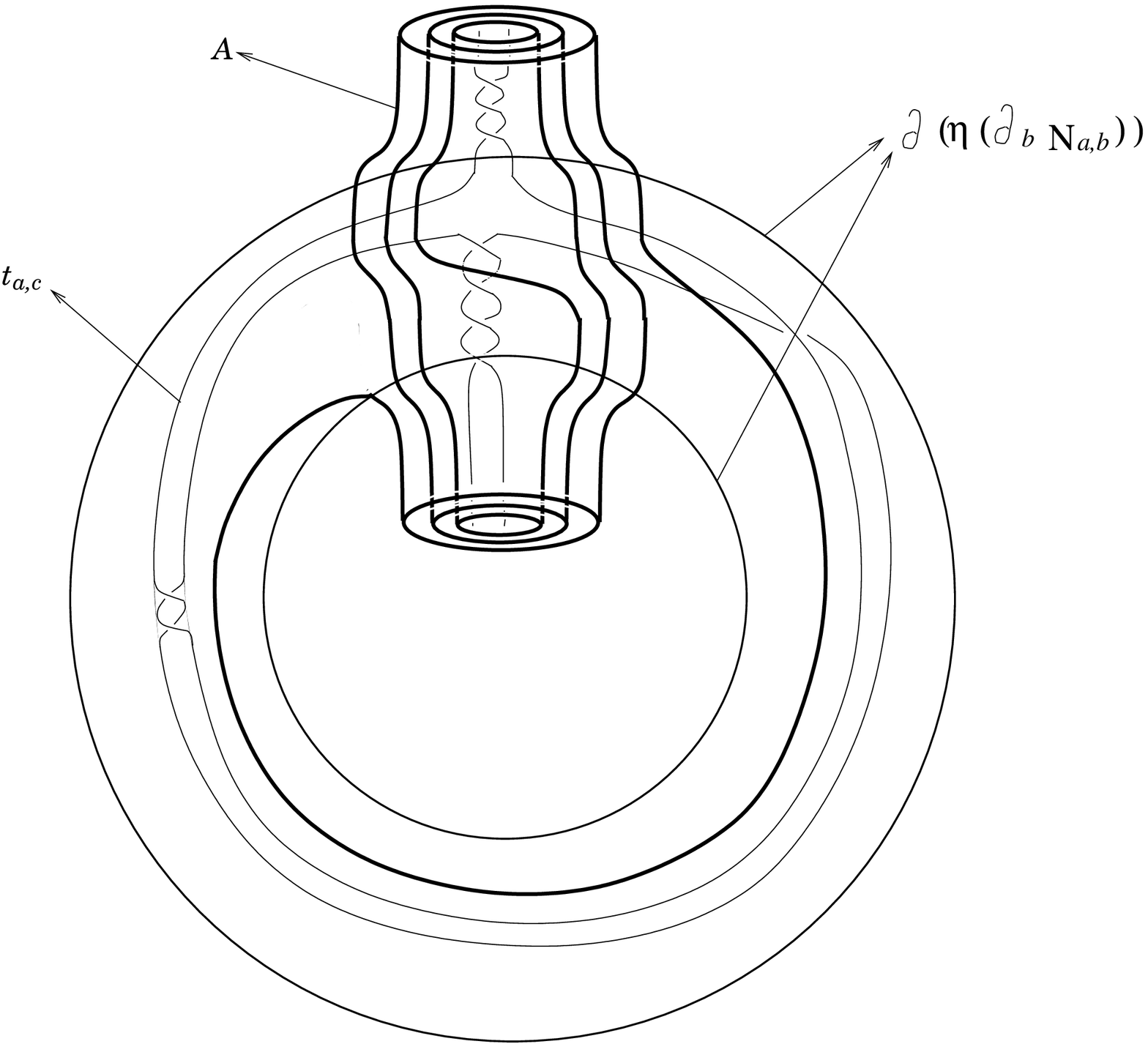}{1}
\endinsert

Suppose that the components of $\partial (\eta(\partial_b N_{a,b}))$, which are
two tori,  lie at levels $x_0$ and $x_{n+1}$, so that
$a < x_0 < x_1< \dots < x_n <x_{n+1} < c$.
Consider the torus $T\times {x_0}$. The annuli $A_i$ intersect this torus
in $2r$ curves which divide it into $2r$ horizontal annuli. Consider the
collection of these annuli, except the one which intersects $t_{a,c}$, and
denote the collection by $B$. Take also an horizontal annulus at a level
between $x_1$ and $x_2$, which does intersect $t_{a,c}$ in one point, and
denote it by $B'$. But if different components of $t_{a,c}$ define the points
$x_1$ and $x_2$, i.e., a component of $t_{a,c}$ goes from $D_1$ to $x_1$, and
the other goes from $D_1$ to $x_2$, then choose $B'$ as an annulus between
$x_2$ and $x_3$; in this case $B'$ is disjoint from $t_{a,c}$. Note that in
this case the points $x_1$ and $x_2$ lie on different vertical annuli,
for otherwise it contradicts condition 2.2.1(2).
Analogously define horizontal annuli
$C$ at level $x_{n+1}$ and an annulus $C'$ at a level between $x_{n}$ and
$x_{n-1}$ (or between $x_{n-1}$ and $x_{n-2}$). If $t_{a,c}$ is disjoint from
$\Cal A'$, then consider only the collection of annuli $B$.

\proclaim{Lemma 2.2.2} The surface $\Cal A$ is incompressible and meridionally
incompressible in $N_{a,c}-int\, \eta (t_{a,c})$.
\endproclaim

\demo{Proof} Suppose first that the arcs $t_{a,c}$ do intersect the surface
$\Cal A'$. Suppose
$\Cal A$ is compressible or meridionally compressible, and suppose that $D$ is a
compression  disk for
$\Cal A$, which is  disjoint from $t_{a,c}$, or intersects it in one point.
Assume $D$ intersects transversely the annuli $B, B', C, C'$.
Let $\gamma$ be a simple closed curve of intersection between $D$ and
the annuli which is innermost
in $D$, and then it bounds a disk $D'$. Note that $\gamma$ must be trivial
in the corresponding annulus, for the
core of the annulus is an essential curve in $N_{a,c}$, and then $\gamma$
bounds a disk $D''$ in such annulus. If $D''$ intersects the arcs $t_{a,c}$
in two points, then $D'\cup D''$ bounds a 3-ball containing part of the arcs,
and one subarc of $t_{a,c}$ will have endpoints on $D''$, which is not possible.
So the disk $D''$ is disjoint from the arcs $t_{a,c}$ or
intersects them in one point, depending if the innermost disk in $D$
intersects or not the arcs $t_{a,c}$; but in any case, by doing an isotopy this
intersection curve is eliminated.

Suppose then that the intersection consists only of arcs, and let
$\gamma$ be an outermost arc in $D$, which we may assume cuts off a disk $D'$
from $D$ which is disjoint from $t_{a,c}$. Let $\partial D' = \gamma\cup
\alpha$. The arc $\gamma$ lies in one of the annuli, and suppose first it is
not an spanning arc in the corresponding annulus. So $\gamma$ bounds a disk
$E$ in such annulus. If $E$ is contained in one of the annuli $B$ or $C$, then
$E$ and $t_{a,c}$ are disjoint, and then by cutting $D$ with $E$ (or with an
innermost disk lying on $E$), we get another compression disk having fewer
intersections with the annuli.
Suppose then that $E$ lies in $B'$ or in $C'$. If
$E$ is disjoint from $t_{a,c}$, the same argument applies, so assume it
intersects $t_{a,c}$. There are two cases, either $D'$ lies in the region between
$B'$ and $D_1$ (or $C'$ and $D_1'$), or it lies in a region bounded by
$B'$, two vertical annuli and $C'$ or one of the annuli $C$ (or bounded by
$C'$, two vertical annuli and $B'$ or one of the annuli $B$).

Suppose first that $D'$ lies in the region between $B'$ and
$D_1$. Let $\hat A_1$ be the part of $A_1$ lying between $T\times \{a\}$
and the level determined by $B'$. The arc $\alpha$ lies on $\hat A_1$, with
endpoints on the same boundary component of $\hat A_1$. Then $\alpha$ is a
separating arc in $\hat A_1$, and cuts off a subsurface $E'$, which is a disk or
an annulus (depending if $E'$ contains or not $\partial_0 A_1$).
If $E'$ is a disk, then $D' \cup E \cup E'$ is a sphere intersecting
the arcs $t_{a,c}$ in one point, which is impossible (note that in this case
the point $x_1$ cannot lie in $E'$). If $E'$ is an annulus, then
$D' \cup E \cup E' \cup D_1$ bound a 3-ball $P$, so that $\partial P$ intersects
the arcs $t_{a,c}$ in at least 3 points, so it must intersect them in 4 points,
and then the point $x_1$ must lie on $E'$. As the arcs lie in the 3-ball $P$,
it is clear that there is a meridian disk of the solid torus bounded
by $D_1\cup \hat A_1 \cup B'$ which is disjoint from the arcs $t_{a,c}$,
contradicting condition 2.2.1(1).
The case when $D'$ lies in the region between $C'$ and $D_1'$  can be handled
in a similar way.

Suppose now that $D'$ lies in a region bounded by $B'$, two vertical annuli,
and $C'$ or one of the annuli $C$. The arc $\alpha$ cuts off a disk $E''$ from
one of the vertical annuli. The disks $E$, $D'$ and $E''$ form a sphere which
intersects  $t_{a,c}$ in at least two points.  If intersects it in more than
2 points, then there is a subarc of $t_{a,c}$ arc going from $E$ to $E''$,
and at least one sub arc with both endpoints on $E''$.
As the arcs are monotonic, these cannot be tangled,
so any arc with both endpoints on $E''$ must be parallel to $E''$, which
contradicts 2.2.1(3). If the sphere $E \cup D' \cup E''$ intersects
$t_{a,c}$ in two points, then there is an arc of
$t_{a,c}$ which intersects $E''$ at the point $x_j$, say.
We can isotope the arc so that intersects the vertical annuli at a level
just below $B'$; this changes the order of the intersection points $x_i$, so that
the point $x_j$ now becomes $x_2$. Note that $x_1$ and $x_2$ lie on
different vertical annuli because of condition 2.2.1(2).
Now $B'$ is disjoint from $t_{a,c}$,
and by an isotopy we get a disk having fewer intersections with $B'$.

We have shown that the outermost arc $\gamma$ in $D$ cannot bound a disk in one of
the annuli. Suppose then that $\gamma$ is an spanning arc of
the corresponding annulus. Remember that $\partial D'=\gamma\cup \alpha$.
If $\gamma$ lies in $B$ or $C$, then $\alpha$ must be
a spanning arc of $A_r$ or $A_r'$,
and then $\gamma \cup \alpha$ is an essential curve in
$N_{a,b}$, so it cannot bound a disk. If $\gamma$ lies in
$B'$, then $\gamma \cup \alpha$ must be a meridian of the solid torus
$D_1\cup \hat A_1 \cup B'$, but then condition 2.2.1(1) is not satisfied.
Analogously, if $\gamma$ lies on $C'$, condition 2.1.1(4) is not satisfied.

Therefore, if there is a
compression or meridional compression disk $D$ for $\Cal A$, it must be disjoint
from the annuli $B, B', C, C'$. If $D$ lies in a  region bounded by two vertical
annuli and $B$ (or $B'$)  and $C$ (or $C'$), then as the arcs $t_{a,c}$ are straight,
it is not difficult to see that there is a subarc parallel to one of the vertical
annuli, which  contradicts condition 2.2.1(3). If $D$
lies in the solid torus determined by $\hat A_1$, $B'$ and $D_1$, then either this
contradicts condition 2.2.1(1), or $D$ is not a compression or
meridional compression disk. If $D$ lies in any of the remaining regions,
then it is not difficult to see that it cannot be a compression disk.

Suppose now that the arcs $t_{a,c}$ are disjoint from $\Cal A'$. If $D$ is a
compression disk for $\Cal A'$ which is disjoint from the arcs $t_{a,c}$ or
intersects them in one point, a similar argument as above shows that $D$ is
disjoint from the annuli $B$, and then $D$ must lie in the
solid torus bounded by $D_1 \cup A_1 \cup A_1' \cup D_1'$.
This implies that the wrapping number of the knot $k$ defined in 2.2.1(6)
is $0$ or $1$, contradicting that condition.
\qed
\enddemo


\medskip
\subhead {2.3 Pieces of type $\tilde \Cal B$} \endsubhead
\medskip

Let $T$ be a torus, and let $N_{a,b}=T\times [a,b]$. Denote by $A_1,\dots, A_r$ a
collection of disjoint, parallel, once punctured annuli properly embedded in
$N_{a,b}$ as in 2.1. Note that the punctures of these annuli, $\partial_0 A_i$,
are a collection of $r$  trivial, nested curves on $\partial_a N_{a,b}$;
each of them bounds a disk $D_i$, such that $D_i\subset D_{i+1}$,
for $i\in \{1,2,\dots,r-1\}$. The curves
$\partial_1 A_i$, are a collection of $2r$ parallel essential curves on
$\partial_b N_{a,b}$. Let
$A_1',...,A_r'$ be a collection of properly embedded parallel  annuli in a
solid torus $R_b$,  such that the boundaries of the annuli are  a collection
of $2r$ essential curves in $\partial R_b $, each one going at least twice
longitudinally around $R_b$.
Assume that $A_1'$ separates a solid torus $N_2 \subset R_b$,
such that $N_2\cap \partial R_b$ is an annulus parallel to $A_1'$,
and that the interior of $N_2$ does not intersect any $A_i'$.
Identify $\partial_b N_{a,b}$ with $\partial
R_b$, so that the collection of curves $\{ \partial_1 A_i \}$ is identified
with the collection $\{ \partial A_i'\}$.
Note that the boundaries of an annulus in one side may
be identified to curves belonging to two different annuli in the other side.
Let $\Cal B'$ be the union of the annuli $A_i$ and $A_i'$. The surface
$\Cal B'$ is  properly embedded in the solid torus $N_{a,b}\cup R_b$. Note that each
component of
$\Cal B'$ is a punctured torus, and that $\Cal B'$ could be connected. Let
$\eta(\partial_b N_{a,b})$ be a regular neighborhood of $\partial_b N_{a,b}$
in $N_{a,b}\cup R_b$ such that $\eta(\partial_b N_{a,b})\cap R_b$ consists
of a product neighborhood  $\partial R_b \times [b,b+\epsilon]$, and
$\eta (\partial_b N_{a,b})\cap {\Cal B'}$ consists of a collection of
$2r$ parallel vertical annuli. Let $t$ be an arc properly embedded in
$N_{a,b}\cup R_b$, such that: $t\cap \partial_a N_{a,b}$ consists of two
points contained in $D_1$, $t\cap N_{a,b}$ consists of two straight arcs
properly embedded in $N_{a,b}$,
$t\cap R_b$ consists of an arc contained in the product neighborhood
$\partial R_b \times [b,b+\epsilon]$, which has only one minimum in there, i.e.,
it intersects each torus twice, except one, which intersects precisely in a
tangency point. Suppose $t$
intersects $\Cal B'$ in finitely many points, all of which lie in
$\eta(\partial_b N_{a,b})$, and suppose these intersection points lie at
different heights; denote the points and the height at
which they occur by $x_1,x_2,\dots,x_n$, which are ordered
according to their heights.

\medskip
\noindent {\bf 2.3.1} Let $\Cal B = \Cal B'- int\, \eta (t)$.
We call a solid torus together with a surface $\Cal B$ and an arc
$t$, a piece of type $\tilde \Cal B$. We assume that a piece of type
$\tilde \Cal B$ satisfies the following:

\roster

\item The part of $A_1$ up to level $x_1$ satisfy the same condition
as in 2.2.1(1).

\item The same condition as in 2.2.1(2) is satisfied.

\item Suppose a subarc of $t$, which does not contain the minimum of $t$,
intersects a component of $\Cal B'$ in two points, say $x_i$, $x_j$, but does
not intersect any other component between these two points. Then this part
of the arc, say $\beta$, is not parallel to an arc $\alpha$ on $\Cal B'$.
That is, if $D$ is any disk in $N_{a,b}\cup R_b$ with interior disjoint from
$\Cal B'$, such that $\partial D=\alpha \cup \beta$, $\partial D\cap \Cal B'=\alpha$,
then $t$ necessarily intersects $int\, D$. If this condition is not
satisfied then the surface $\Cal B$ will be clearly compressible.

\item Assume that the subarc of $t$ which contains the minimum is
contained in the solid torus $N_2$,
and that the point $x_n$ lies on the curve $A_1'\cap \partial R_b$. The endpoints
of the arc $t\cap R_b$ lie on the annulus $N_2\cap \partial R_b$.
Suppose there is
no a meridian disk $D$ of $N_2$, disjoint from $t$, and whose boundary
consists of an arc on $A_1'$, and one arc in $\partial R_b$.

\noindent Note that if this condition is not satisfied, and so there is such a meridian disk,
then the subarc $t'$ of $t$ lying in $N_2$,
is parallel onto $\partial N_2$, that is, there is a disk $F$ embedded in $N_2$,
with $\partial F = t'\cup \alpha\cup \beta$, where $\alpha$ is an arc on $A_1'$
and $\beta$ an arc on $\partial R_b$, and $int\ F\cap t'=\emptyset$.
So we can slide the arc $t'$ until it is contained in $\partial R_b$;
by sliding it further, we could find a compression disk for $\Cal A$.

\item A condition similar to 2.2.1(4) is satisfied for $x_n$. That is, one subarc of
$t$ starts at the maximum of $t$ and arrives to the point $x_n$. Let
$t'$ be the other subarc of $t$ that starts at the maximum. Assume that either (i)
$t'$ is disjoint from $\Cal B'$; or (ii) $x_n$ and the first point of
intersection of  $t'$ with $\Cal B'$, say $x_j$, lie on different vertical annuli
of $\eta(\partial_b N_{a,b})\cap\Cal B'$; or (iii) $x_n$ and $x_j$ lie on the
same vertical annulus, but there is no disk
$D \subset \eta(\partial_b N_{a,b})$, with
$\partial D =\alpha \cup \beta\cup\gamma$, so that $\alpha\subset t'$,
$\beta \subset \Cal B'$, which connects $x_j$ with a point $x_n'$ at level
$T\times \{x_n\}$, $\gamma$ is an arc on
$T\times \{x_n\}$, and $int\, D$ is disjoint from $t$ and from $\Cal B'$. In
other words, $t'$ cannot be slided to lie at level $T\times \{x_n\}$.

\noindent This condition complements condition (4), for if this is not satisfied,
the arcs could be slided into a position in which (4) fails.

\item In case $t$ is disjoint from $\Cal B'$, we assume the following:
In this case the annulus $A_i$ is glued to the annulus $A_i'$, for all $i$. Let $k$ be
the knot obtained by joining the endpoints of $t$ lying on $D_1$ with an arc
on $D_1$. Then we get a knot lying in the solid torus bounded by $A_1 \cup
A_1'$. Assume that this knot has wrapping number $\geq 2$ in such a solid
torus. If this is not satisfied, then the surface will be compressible or
meridionally compressible.
\endroster

It is not difficult to construct
pieces of type $\tilde \Cal B$ where all of these conditions are satisfied.
See Figure 2.

If the annulus $A_1$ goes just one longitudinally around the solid torus $R_b$,
and $t$ intersects $\Cal B'$, then it is not difficult to see that $\Cal B$
will be compressible; if $t$ is disjoint from $\Cal B'$, then $\Cal B$ is
incompressible but each of its components will be parallel into
$\partial_a N_{a,b}$.

\remark {Remark} Condition 2.2.1(4) has an alternative description.
Embed the solid torus $N_2$
in $S^3$ in a standard manner, so that $A_1'\cap \partial R_b$
is a preferred longitude of such solid torus. Suppose the point $x_n$ lies
on the curve $A_1'\cap \partial R_b$. As before, the endpoints of the arc
$t\cap R_b$ lie on the annulus $N_2\cap \partial R_b$  and
the point $x_n$ lies on the curve $A_1'\cap \partial R_b$
Let $\alpha$ be the boundary of a meridian disk of $N_2$ which passes through
these points.
The endpoints of $t\cap R_b$ separates $\alpha$ into two arcs;
connect them with the subarc of $\alpha$ which intersects $A_1'$.
This defines a knot $k$
in $S^3$, which by construction has a presentation with two minima.
It can be shown that Condition 2.2.1(4) is satisfied if and only if
$k$ is a non-trivial 2-bridge knot. Note that any 2-bridge knot $k$
can result from this construction, as shown in Figure 3.

\midinsert
\figure{3}{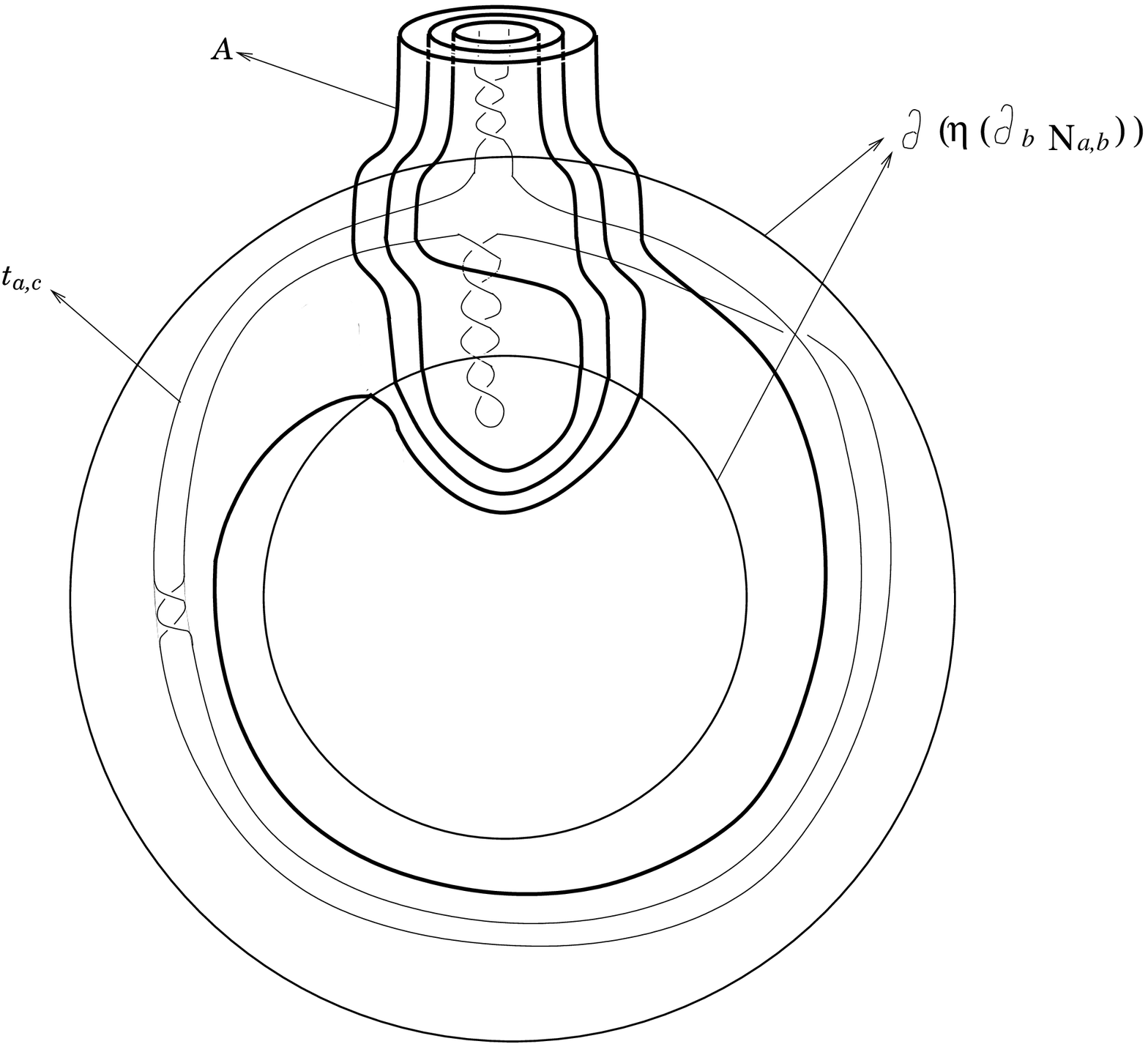}{2}
\endinsert

\midinsert
\figure{2.5}{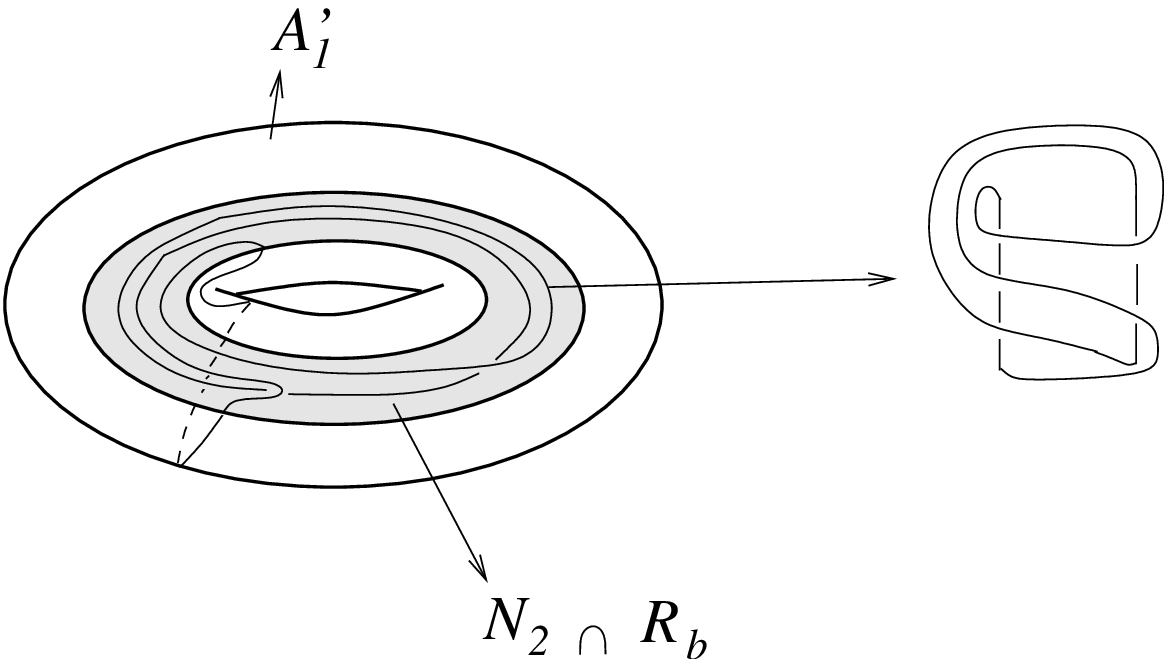}{3}
\endinsert

Suppose that the components of $\partial (\eta(\partial_b N_{a,b}))$ lie at
levels $x_0$ and $x_{n+1}$, so that
$a < x_0 < x_1<\dots x_n <x_{n+1}\leq b+\epsilon$.
Consider the torus $T\times \{x_0\}$; note that ${\Cal B} \cap (T\times \{x_0\})$
consist of $2r$ curves which divide that torus into $2r$ horizontal annuli.
Consider the
collection of these annuli, except the one which intersects $t$, and
denote the collection by $B$. Take also an horizontal annulus at a level
between $x_1$ and $x_2$, which intersects $t$ in one point, and
denote it by $B'$.
But if different subarcs of $t$ define the points
$x_1$ and $x_2$, i.e., a subarc of $t$ goes from $D_1$ to $x_1$, and
the other goes from $D_1$ to $x_2$, then choose $B'$ as an annulus between
$x_2$ and $x_3$; in this case $B'$ is disjoint from $t$, and the points
$x_1$, $x_2$ lie on different vertical annuli.
Analogously define horizontal annuli
$C$ at level $x_{n+1}$ and annulus $C'$ at a level between $x_{n}$ and
$x_{n-1}$ (or between $x_{n-1}$ and $x_{n-2}$) . If $t$ is disjoint from
$\Cal B'$, then consider only the collection of annuli $B$.

\proclaim{Lemma 2.3.2} The surface $\Cal B$ is incompressible and meridionally
incompressible in $(N_{a,b}\cup R_b)- int\, \eta (t)$.
\endproclaim

\demo{Proof} Suppose first that the arc $t$ does intersect the surface $\Cal B'$.
Suppose $D$ is a compression  disk for $\Cal B$, which is
disjoint from $t$, or intersects it in one point. Suppose $D$
intersects transversely the annuli $B, B', C, C'$. Let $\gamma$ be a simple
closed curve of intersection between $D$ and the annuli which is innermost
in $D$. Note that $\gamma$ must be trivial in the corresponding annulus, for
the core of the annulus is an essential curve in $T\times [a,b]$ or in $R_b$.
So $\gamma$ bounds a disk in such annulus, which is disjoint from the
arc $t$ or intersects it in one point, depending if the innermost disk in
$D$ intersects or not the arc $t$, but in any case an isotopy eliminates this
intersection.

Suppose then that the intersection between $D$ and the annuli consists
only of arcs, and let
$\gamma$ be an outermost arc in $D$, which we may assume it cuts off a disk $D'$
from $D$ which is disjoint from $t$. Let $\partial D' = \gamma\cup \alpha$.
The arc $\gamma$ lies in one of the annuli, and suppose first it is
not an spanning arc in the corresponding annulus. So $\gamma$ bounds a disk
$E$ in such annulus. If $E$ is contained in one of the annuli $B$ or $C$, then
$E$ and $t$ are disjoint, and then by cutting $D$ with $E$ (or with an
innermost disk lying on $E$), we get another compression disk having fewer
intersections with the annuli.
Suppose then that $E$ lies in $B'$ or in $C'$. If
$E$ is disjoint from $t$, the same argument as above applies, so assume it
intersects $t$. There are three cases, either $D'$ lies in the region between
$B'$ and $D_1$; or it lies in a region bounded by
$B'$, two vertical annuli and $C'$ or one of the annuli $C$; or it lies
in the region bounded by $C'$ and $A_1'$.

Suppose first that $D'$ lies in the region between $B'$ and
$D_1$. Let $\hat A_1$ be the part of $A_1$ lying between $T\times \{a\}$
and the level determined by $B'$. The arc $\alpha$ lies on $\hat A_1$, with
endpoints on the same boundary component of $\hat A_1$. Then $\alpha$ is a
separating arc in $\hat A_1$, and cuts off a subsurface $E'$, which is a disk or
an annulus (depending if $E'$ contains or not $\partial_0 A_1$).
If $E'$ is a disk, then $D' \cup E \cup E'$ is a sphere intersecting
the arcs $t$ in one point, which is impossible (note that in this case
the point $x_1$ cannot lie in $E'$). If $E'$ is an annulus, then
$D' \cup E \cup E' \cup D_1$ bound a 3-ball $P$, so that $\partial P$ intersects
the arcs $t$ in at least 3 points, so it must intersect them in 4 points,
and then the point $x_1$ must lie on $E'$.
As the arcs lie in the 3-ball $P$,
it is clear that there is a meridian disk of the solid torus bounded
by $D_1\cup \hat A \cup B'$ which is disjoint from the arcs $t$,
contradicting condition 2.3.1(1)

Suppose now that $D'$ lies in a region bounded by $B'$, two vertical annuli,
and $C'$ or one of the annuli $C$. The arc $\alpha$ cuts off a disk $E''$ from
one of the vertical annuli. The disks $E$, $D'$ and $E''$ form a sphere which
intersects  $t$ in at least two points.  If intersects it in more than
2 points, then there is one subarc of $t$ going from $E$ to $E''$, and at least one
arc with both endpoints on $E''$. As the arcs are monotonic, these cannot be tangled,
so any arc with both endpoints on $E''$ must be parallel to $E''$, which
contradicts 2.3.1(2).
If the sphere $E \cup D' \cup E''$ intersects
$t$ in two points, then there is an arc of
$t$ which intersects $E''$ at the point $x_j$, say.
We can isotope the arc so that intersects the vertical annuli at a level
just below $B'$; this changes the order of the intersection points $x_i$, so that
the point $x_j$ now becomes $x_2$. Note that $x_1$ and $x_2$ lie on
different vertical annuli because of condition 2.3.1(2).
Now $B'$ is disjoint from $t$,
and by an isotopy we get a disk having fewer intersections with $B'$.

Suppose then that $D'$ lies in the region bounded by $C'$ and $A_1'$.
The arc $\alpha$ cuts off a disk $E'$ from $A_1'$. Note that $D'$
does not intersect the arc $t$, so $E$ and $E'$ both intersect $t$
in one point, and $D'\cup E \cup E'$
bounds a 3-ball containing a subarc of $t$. Then there is a meridian
disk of the solid torus bounded by $C'$ and $A_1'$ disjoint from $t$,
but this contradicts 2.3.1(4).

We have shown that the outermost arc $\gamma$ in $D$ cannot bound a disk in one of
the annuli. Suppose then that $\gamma$ is a spanning arc of
the corresponding annulus. Remember that $\partial D'=\gamma\cup \alpha$.
If $\gamma$ lies in $B$ or $C$, then $\alpha$ must be a spanning arc of $A_r$
or $A_r'$,and then $\gamma \cup \alpha$ is an essential curve in
$N_{a,b}$, so it cannot bound a disk (the curve $\gamma \cup \alpha$
could bound a disk only if $A_r'$ were a longitudinal annulus in $R_b$).
If $\gamma$ lies in $B'$, then $\gamma \cup \alpha$ must be a meridian
of the solid torus bounded by $D_1\cup \hat A_1 \cup B'$, but then condition
2.3.1(1) is not satisfied.
Analogously, if $\gamma$ lies on $C'$, condition 2.3.1(4) is not satisfied.

Therefore, if there is a compression or
meridional compression disk for $\Cal B$, it must be disjoint from the annuli
$B, B', C, C'$. If $D$ lies in a  region bounded by two vertical
annuli and $B$ (or $B'$)  and $C$ (or $C'$), then as the arcs $t$ are straight,
it is not difficult to see that there is a subarc parallel to one of the vertical
annuli, which  contradicts condition 2.3.1(3). If $D$
lies in the solid torus determined by $\hat A_1$, $B'$ and $D_1$, then either this
contradicts condition 2.3.1(1), or $D$ is not a compression or
meridional compression disk. If $D$ lies in any of the remaining regions,
then it is not difficult to see that it cannot be a compression disk.

The remaining case is when $t$ is disjoint from $\Cal B'$. In this
case the annulus $B'$ is not defined. Let $D$ a compression or
meridional compression disk. Again, we can make $D$ disjoint from
the annuli $B$. Then $D$ lies in the solid torus determined by
$A_1\cup A_1'$. If $\partial D$ is essential in that torus, this
implies that the wrapping number of the knot $k$ defined in 2.3.1(6)
is $0$ or $1$, contradicting that condition. \qed
\enddemo

\medskip
\subhead {2.4 Pieces of type $\tilde \Cal C$} \endsubhead
\medskip

Let $T$ be a torus, and let $N_{a,b}=T\times [a,b]$. Denote by
$A_1,\dots, A_r$ a collection of disjoint, parallel, once
punctured annuli properly embedded in $N_{a,b}$ as in 2.1. Note that the
punctures of these annuli, $\partial_0 A_i$, are a collection of $r$
trivial, nested curves in $\partial_a N$, each bounding a disk $D_i$,
and $\partial_1 A_i$ consists of a
collection of $2r$ essential parallel curves on $\partial_b N_{a,b}$. Let
$A_1',...,A_{2r}'$ be a collection of meridian disks in a solid torus $R_b$.
Identify $\partial_b N_{a,b}$ with $\partial R_b$, so that the collection
of curves $\{\partial_1 A_i\}$ are identified with the collection
$\{ \partial A_i'\}$,
and say, the annulus $A_i$ is glued with the disks $A_i'$ and $A_{2r-i+1}'$.
Let $\Cal C'$
be the union of the annuli $A_i$ and the disks $A_i'$. $\Cal C'$ is a surface
properly embedded in the solid torus $N_{a,b}\cup R_b$.
Note that each component of $\Cal C'$
is a disk.  Let $\eta(\partial_b N_{a,b})$ be a regular neighborhood of
$\partial_b N_{a,b}$ in $N_{a,b}\cup R_b$ such that
$\eta(\partial_b N_{a,b})\cap R_b$ consists
of a product neighborhood  $\partial R_b \times [b,b+\epsilon]$, and
$\eta (\partial_b N_{a,b})\cap {\Cal C'}$ consists of a collection of
$2r$ parallel vertical annuli. Let $t$ be an arc properly embedded in
$N_{a,b}\cup R_b$, such that: $t\cap \partial_a N_{a,b}$ consists of two
points contained in the disk $D_1$, $t\cap N_{a,b}$ consists of two straight arcs
properly embedded in $N_{a,b}$, $t\cap R_b$ consists of an arc contained
in the product neighborhood
$\partial R_b \times [b,b+\epsilon]$, which has only one minimum in there,
i.e., it intersects each torus twice, except one, which intersects precisely in
a tangency point.  Assume that the arc $t$ intersects
$\Cal C'$ in finitely many points, all of which lie in $\eta(\partial_b
N_{a,b})$. Suppose these intersection points lie at different heights,
are denoted by $x_1,x_2,\dots,x_n$, and are ordered according to their heights.

\medskip
\noindent {\bf 2.4.1} Let $\Cal C = \Cal C'- int\, \eta (t)$.
We call a solid torus together with a surface $\Cal C$ and an arc
$t$, a piece of type $\tilde \Cal C$. We assume that a piece of type
$\tilde \Cal C$ satisfies the following:

\roster

\item The part of $A_1$ up to level $x_1$ satisfies the same condition as
in 2.2.1(1).

\item A similar condition as in 2.2.1(2) is satisfied.

\item Suppose a subarc of $t$, say $\beta$, which does not contain the
minimum of $t$,
intersects a component of $\Cal C'$, in two points, say $x_i$, $x_j$, but
does not intersects any other component between these two points. The same
condition as in 2.3.1(3) is satisfied, that is, if $D$ is any disk in
$N_{a,b}\cup R_b$ with interior
disjoint from $\Cal C'$, such that $\partial D=\alpha \cup \beta$,
$\partial D\cap \Cal C'= \alpha$, then $t$
necessarily intersects $int\, D$.

\item The subarc of $t$ which contains the minimum has endpoints on
two different disks $A_i'$ and $A_j'$.
\endroster

Note that in this case the arc $t$ necessarily intersects
the surface $\Cal C'$.
It is not difficult to construct examples
of pieces satisfying these conditions. See Figure 4.

\midinsert
\figure{3}{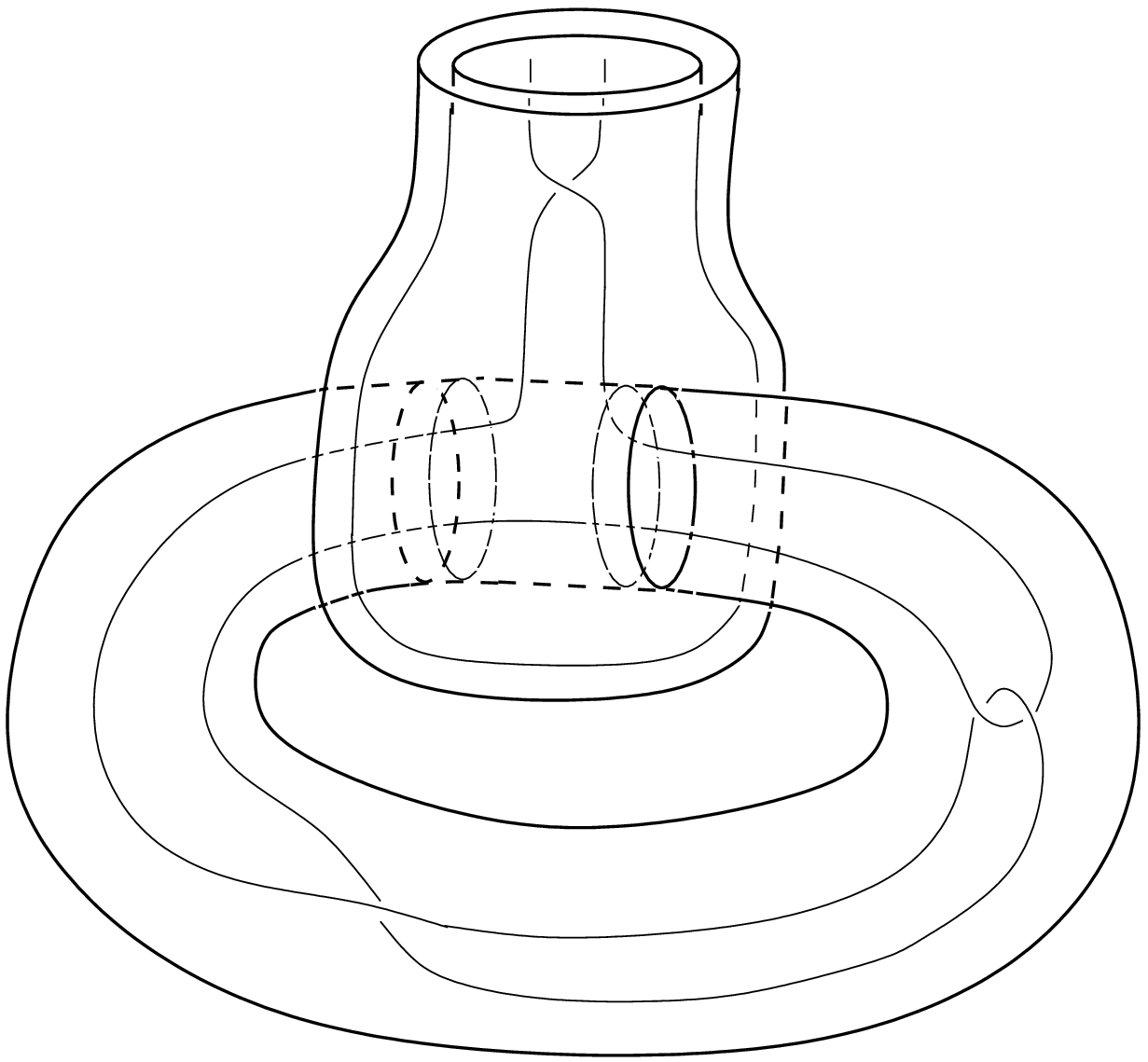}{4}
\endinsert

Suppose that the components of $\partial (\eta(\partial_b N_{a,b}))$ lie at
levels $x_0$ and $x_{n+1}$, so that
$a < x_0 < x_1<\dots x_n <x_{n+1}\leq b+\epsilon$.
Consider the torus $T\times \{x_0\}$, and note that ${\Cal C} \cap (T\times \{x_0\})$
consists of $2r$ curves which divide that torus into $2r$ horizontal annuli.
Consider the
collection of these annuli, except the one which intersects $t$, and
denote the collection by $B$. Note that in this case the core of any of the annuli $B$
is a trivial curve in $N_{a,b}\cup R_b$. Take also an horizontal annulus at a level
between $x_1$ and $x_2$, which intersects $t$ in one point, and
denote it by $B'$. But if different subarcs of $t$ define the points
$x_1$ and $x_2$, i.e., a subarc of $t$ goes from $D_1$ to $x_1$, and
the other goes from $D_1$ to $x_2$, then choose $B'$ as an annulus between
$x_2$ and $x_3$; in this case $B'$ is disjoint from $t$, and the points
$x_1$, $x_2$ lie on different vertical annuli.

\proclaim{Lemma 2.4.2} The surface $\Cal C$ is incompressible and meridionally
incompressible in $N_{a,b}\cup R_b -int\, \eta (t)$.
\endproclaim

\demo{Proof} Suppose $D$ is a compression  disk for $\Cal C$, which is
disjoint from $t$, or intersects it in one point.

Suppose first that the disk is not in the 3-ball bounded by
$A_1\cup A_1'\cup A_{2r}'$.
Suppose $D$
intersects transversely the annuli $B$. Let $\gamma$ be a simple
closed curve of intersection between $D$ and the annuli which is innermost
in $D$. If $\gamma$ bounds a disk in one of the annuli, then this disk
is disjoint from the arc $t$ and an isotopy eliminates this
intersection. Suppose then that $\gamma$ is an essential curve in one of the
annuli $B$; if this happens then the arc $t$ is disjoint from a meridian
of the solid torus $N_{a,b}\cup R_b$ or it intersects a meridian in one point.
This implies that condition 2.4.1(3) is not satisfied, except if $t$ intersects
each of the disks $A_i\cup A_i' \cup A_{2r-i+1}'$ in two points, but in this
case it is not difficult to see
that $\partial D$ would not be essential in $\Cal C'$.

Suppose then that the intersection between the disk $D$ and the annuli $B$
consists only of arcs, and let
$\gamma$ be an outermost arc in $D$, which we may assume it cuts off a disk $D'$
from $D$ which is disjoint from the arc $t$. Let
$\partial D' = \gamma\cup\alpha$. The arc $\gamma$ lies in one of the
annuli $B$, and suppose first it is not an spanning arc in the corresponding
annulus. So $\gamma$ bounds a disk
$E$ in such annulus. As $E$ is contained in an annulus in $B$,
$E$ is disjoint from $t$, and then by cutting $D$ with $E$ (or with an
innermost disk lying on $E$), we get another compression disk having fewer
intersections with the annuli. Suppose then that $\gamma$ is an spanning arc of
the corresponding annulus. Then $\alpha$ must be an arc contained in one of the
components of $\Cal C$, and the only possibility is that it is in the
component formed by $A_r\cup A_r' \cup A_{r+1}'$, for it is the only component
which intersects a single annulus of $B$ in two curves. But then $\gamma\cup \alpha$
is an essential curve in $N_{a,b}$, which cannot bound a disk.
If the intersection between the disk $D$ and the annuli $B$ is empty,
then it is not difficult to see either condition 2.4.1(3) is not satisfied,
or $\partial D$ is not essential in $\Cal C'$.

Suppose now that the disk $D$ is contained in the ball bounded by
$A_1\cup A_1'\cup A_{2r}'$. Assume $D$
intersects transversely the annulus $B'$. Note that the annulus $B'$ plus two disks
contained in $A_1\cup A_1' \cup A_{2r}'$ bound a 3-ball $B''$ which contains all the
points of intersection between $t$ and $A_1\cup A_1' \cup A_{2r}'$, except $x_1$
(and $x_2$, in case $B'$ is disjoint from $t$).
Suppose there is an outermost arc of intersection $\gamma$ in $D$,
which cuts off a disk $D'$ from $D$, with interior disjoint from $B'$ and $t$.
Suppose that $\gamma$ is not an spanning arc in $B'$; so it cuts off a disk $E$ from $B'$.
If $E$ is disjoint from $t$, by cutting $D$ with $E$, we get a compression disk
having fewer intersections with $B'$. So suppose that $E$ intersects $t$ in one point.

There are two cases, either $D'$ lies in the region between $B'$ and $D_1$,
or it lies in the 3-ball $B''$.
Suppose first that $D'$ lies in the region between $B'$ and
$D_1$. Let $\hat A_1$ be the part of $A_1$ lying between $T\times \{a\}$
and the level determined by $B'$. The arc $\alpha$ lies on $\hat A_1$, with
endpoints on the same boundary component of $\hat A_1$. Then $\alpha$ is a
separating arc in $\hat A_1$, and cuts off a subsurface $E'$, which is a disk or
an annulus (depending if $E'$ contains or not $\partial_0 A_1$).
If $E'$ is a disk, then $D' \cup E \cup E'$ is a sphere intersecting
the arcs $t$ in one point, which is impossible (note that in this case
the point $x_1$ cannot lie in $E'$). If $E'$ is an annulus, then
$D' \cup E \cup E' \cup D_1$ bound a 3-ball $P$, so that $\partial P$ intersects
the arcs $t$ in at least 3 points, so it must intersect them in 4 points,
and then the point $x_1$ must lie on $E'$.
As the arcs lie in the 3-ball $P$,
it is clear that there is a meridian disk of the solid tours bounded
by $D_1\cup \hat A \cup B'$ which is disjoint from the arcs $t$,
contradicting condition 2.4.1(1)

Suppose now that $D'$ lies in the 3-ball $B''$. The arc $\alpha$ cuts off a disk $E''$ from
one of the vertical disks. The disks $E$, $D'$ and $E''$ form a sphere which
intersects  $t$ in at least two points.  If intersects it in more than
2 points, then there is a subarc of $t$ going from $E$ to $E''$, and at least one
arc with both endpoints on $E''$. As the arcs are monotonic, these cannot be tangled,
so any arc with both endpoints on $E''$ must be parallel to $E''$, which
contradicts 2.4.1(3). If the sphere $E \cup D' \cup E''$ intersects
$t$ in two points, then there is a subarc of
$t$ which intersects $E''$ at the point $x_j$, say.
We can isotope the arc so that intersects the vertical disk at a level
just below $B'$; this changes the order of the intersection points $x_i$, so that
the point $x_j$ now becomes $x_2$.
The points $x_1$ and $x_2$ lie on
different vertical annuli because of condition 2.4.1(2).
Now $B'$ is disjoint from $t$,
and by an isotopy we get a disk having fewer intersections with $B'$.

We have shown that the outermost arc $\gamma$ in $D$ cannot bound a disk in $B'$.
Suppose then that $\gamma$ is an spanning arc of $B'$.
The only possibility is that $D'$ lies in the region between $B'$ and $D_1$.
Then $D'$ must be a meridian of the solid torus bounded by
$D_1\cup \hat A_1\cup B'$, but then condition 2.4.1(1) is not satisfied.

We have shown that there are no outermost arcs of intersection in $D$,
so the intersection between $D$ and $B'$ is either empty or contains simple closed curves
(and possibly some arcs).
Let $\gamma$ be a simple closed curve of intersection between $D$ and $B'$ which is innermost
in $D$, and bounds a disk $D'\subset D$. If $\gamma$ bounds a disk in $B'$, then this disk
is disjoint from the arc $t$  or intersects it in one point,
depending if
$D'$ intersects or not the arc $t$, but in any case an isotopy eliminates this
intersection. If $\gamma$ is essential in the
annulus $B'$, then $D'$ must be contained in the 3-ball $B''$, and then
condition 2.4.1(3) is not satisfied, unless the arc $t$ intersects
$A_1\cup A_1'\cup A_{2r}'$ in 2 or 4 points. If there are two points of intersection,
just $x_1$ and $x_2$, then $\partial D$ is either parallel to
$\partial (A_1\cup A_1'\cup A_{2r}')$, which implies that the arc $t$ is disjoint from
$A_1\cup A_1'\cup A_{2r}'$, or $\partial D$ encloses a disk containing just one of the points of $t$,
which implies that $\partial D$ is not essential.
Suppose then that $t$ intersects  $A_1\cup A_1'\cup A_{2r}'$ in 4 points;
in this case $D'$ intersects $t$. Any other closed curve of intersection
between $D$ and $B'$ must be concentric with $\gamma$ in $D$, for otherwise there is an
innermost curve bounding a disk disjoint from $t$. This implies that there is no arc of intersection
between $D$ and $B'$, for there will be an outermost one. So the intersection
consists only of curves. As these are concentric, if there is more than one
these can be removed by an isotopy. So, the intersection between $D$ and $B'$
consists
only of the curve $\gamma$. The only possibility is that $\partial D$ is in fact
parallel to $\partial B'$, and then the disk is parallel to  a disk lying
on $A_1\cup A_1'\cup A_{2r}'$, a contradiction.

The only possibility left is that $D$ is disjoint from $B'$. Then it is not difficult
to see that $D$ is parallel to a disk in $A_1\cup A_1'\cup A_{2r}'$, a contradiction.
\qed
\enddemo


\medskip
\subhead {2.5 Pieces of type $\tilde \Cal D$} \endsubhead
\medskip

Let $R_a$ and $R_b$ be two solid tori.
Let $A_1,...,A_r$ be a collection of properly embedded parallel  annuli in
the solid torus $R_a$,  such that the boundaries of the annuli are  a
collection of $2r$ curves in $\partial R_a$  which go at least twice
longitudinally around $R_a$. Let $A_1',...,A_r'$ be a collection of properly
embedded parallel annuli in the solid torus $R_b$, such that the boundaries
of the annuli are a collection of $2r$ curves in $\partial R_b $ which go
at least twice longitudinally around $R_b$. Identify $\partial R_a$ with
$\partial R_b$, so that the collection of curves $\{\partial A_i\}$ are
identified with the collection of curves $\{ \partial A_i'\}$.
Note that the boundaries of an annulus in one
side may be identified to curves belonging to two different annuli. Let
$\Cal D'$ be the union of the annuli $A_i$ and $A_i'$. $\Cal D'$ is a surface
embedded in $R_a \cup R_b$. Note that each component of $\Cal D'$ is
a torus and that $\Cal D'$ could be connected. Let $\eta(\partial R_a)$ be a
regular neighborhood of $\partial R_a$ in $R_a\cup R_b$ such that
$\eta (\partial R_a) \cap R_a$ is a product neighborhood
$\partial R_a\times [-\epsilon,0]$, and  $\eta (\partial R_a) \cap R_b$ is a product
neighborhood $\partial R_a\times [0,\epsilon]$. Assume that
$\eta (\partial R_a) \cap {\Cal D'}$ consists of a collection of
$2r$ parallel vertical annuli. Let $t$ be a knot  embedded in $R_a \cup R_b$,
such that: $t\cap R_b$ consists of an arc properly embedded in the product
neighborhood $\partial R_a\times [0,\epsilon]$ containing
only a minimum (as defined as in 2.3),  $t\cap R_a$ consists of an arc properly
embedded in the product neighborhood $\partial R_a\times [-\epsilon,0]$,
containing only a maximum (defined analogously to a minimum),
and $t$ intersects $\Cal D'$ in finitely many points, all contained in
$\eta (\partial R_a) \cap {\Cal D'}$. Suppose these intersection
points lie at different heights (in $\partial R_a\times [-\epsilon,\epsilon]$),
are denoted  by $x_1,x_2,\dots,x_n$, and are ordered according to their heights.

\medskip
\noindent {\bf 2.5.1} Let $\Cal D = \Cal D'- int\, \eta (t)$.
We call a solid torus together with a surface $\Cal D$ and an arc
$t$, a piece of type $\tilde \Cal D$. We assume that a piece of type
$\tilde \Cal D$ satisfies the following:

\roster

\item Note that $A_1$ separates a solid torus $N_1 \subset R_a$. Assume that
the subarc of $t$ which contains the maximum is contained in $N_1$. A similar
condition as in 2.3.1(4) is satisfied.

\item For $A_1$ and the minimum of $t$, a similar condition as in
2.3.1(5) is satisfied.

\item Suppose an arc, which does not contain the minimum or maximum of $t$,
intersects a component of $\Cal D'$  in two points, say $x_i$, $x_j$, but
does not intersects any other component between these two points. The same
condition as in 2.3.1(3) is satisfied.

\item The annulus $A_1'$ separates a solid torus $N_2 \subset R_b$.
Assume that the subarc of $t$ which contains the minimum is contained in $N_2$.
The same condition as in 2.3.1(4) is satisfied.

\item For $A_1'$ and the maximum of $t$, a similar condition as in
2.3.1(5) is satisfied.

\item If $t$ is disjoint from $\Cal D'$, then  $A_i$ is glued to the annulus
$A_i'$ and $t$ is a knot lying in the solid torus bounded by
$A_i \cup  A_i'$. In this case assume that the knot $t$
has wrapping number $\geq 2 $ in such a solid torus.

\endroster

Note that if one of the annuli $A_i$ or $A_i'$ is longitudinal in $R_a$
or $R_b$, then the surface $\Cal D$ will be compressible.

The space $R_a\cup R_b$ can be $S^3$, $S^1\times S^2$ or a lens
space $L(p,q)$. It is not difficult to construct examples of knots
and surfaces satisfying these conditions for any of the possible
spaces $R_a\cup R_b$. See Figure 5 for a specific example of such a
piece in $S^3$,  where $\Cal D'$ is made of two annuli, i.e., $\Cal
D' = A_1 \cup A_1'$ (shown with gray lines), and say, $R_b$ is the
solid torus standardly embedded in $S^3$, and $R_a$ is the
complementary torus, which contains the point at infinity. The knot
$t$ intersects $\Cal D'$ in two points, which divide it into two
arcs; the black bold arc is the one which contains the maximum of
$t$. Note that in this example the torus $\Cal D'$ is in fact
isotopic to $\partial R_a$ in $S^3$, but $\Cal D$ is incompressible
and meridionally incompressible in $S^3-int\, \eta (t)$, while
$\partial R_a-int\, \eta (t)$ is compressible.

\midinsert
\figure{2.9}{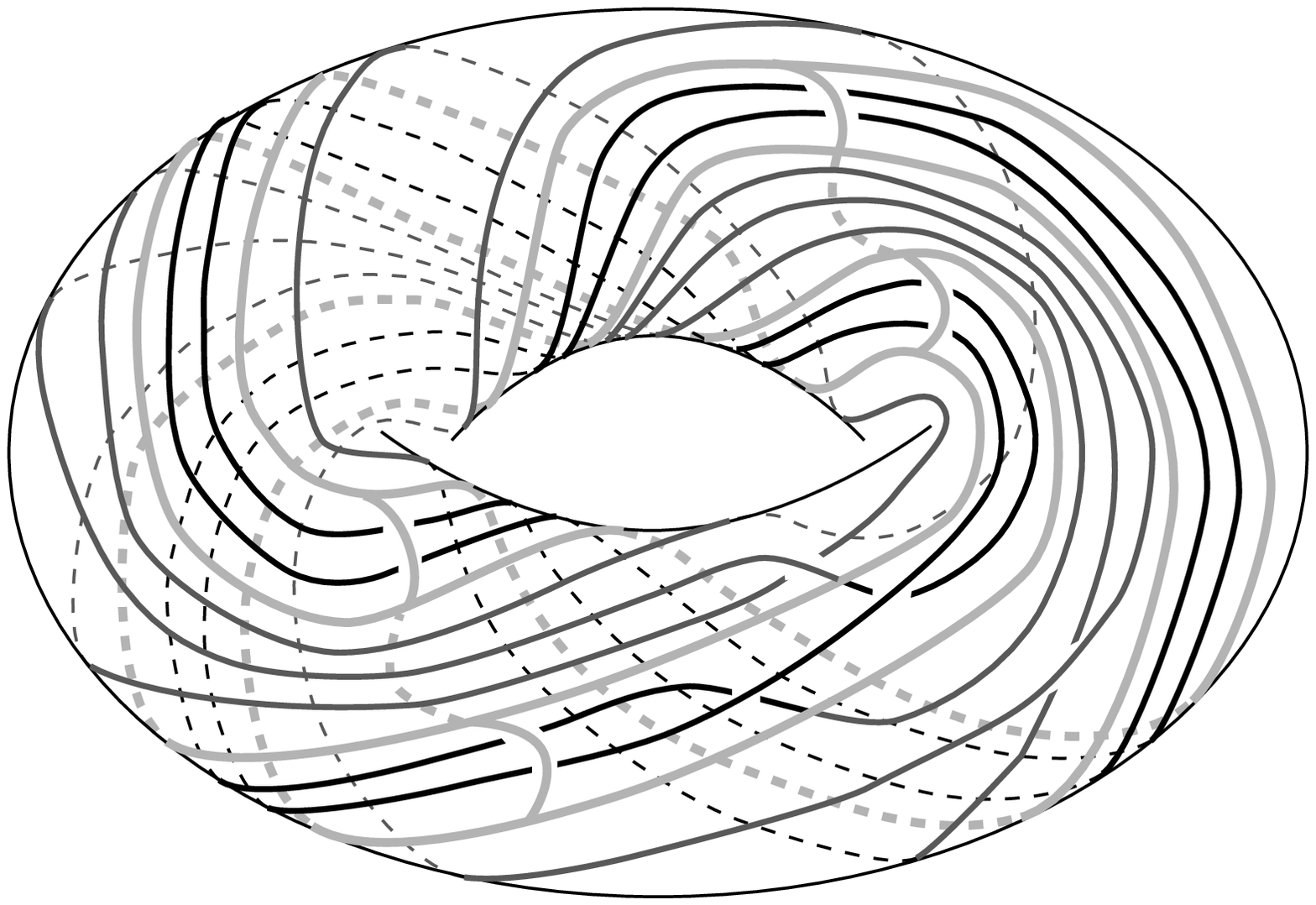}{5}
\endinsert

\proclaim{Lemma 2.5.2} The surface $\Cal D$ is incompressible and meridionally
incompressible in $R_a \cup R_b -int\, \eta (t)$.
\endproclaim

\demo{Proof} It is similar to the proof of Lemma 2.3.2.
\enddemo

\medskip
\subhead {2.6  Pieces of type $\tilde \Cal E$} \endsubhead
\medskip

Let $R_a$ and $R_b$ be two solid tori.
Let $A_1,...,A_r$ be a collection of properly embedded parallel annuli in
the solid torus $R_a$, such that the boundaries of the annuli are  a
collection of $2r$ curves in $\partial R_a$ which go at least twice
longitudinally around $R_a$. Let $A_1',...,A_{2r}'$ be a collection of
meridian disks in a solid torus $R_b$. Identify $\partial R_a$ with
$\partial R_b$, so that the collection of curves $\{\partial A_i\}$
are identified with the collection $\{\partial A_i'\}$.
Let $\Cal E'$ be the union of the annuli $A_i$
and the disks $A_i'$. $\Cal E'$ is a surface embedded in $R_a\cup R_b$ and
each of its components is a sphere. Note that $R_a\cup R_b$ can be
any lens space $L(p,q)$, but it cannot be $S^3$ nor $S^1\times S^2$.
Let $\eta(\partial R_a)$ be a
regular neighborhood of $\partial R_a$ in $R_a\cup R_b$ such that
$\eta (\partial R_a) \cap R_a$ is a product neighborhood
$\partial R_a\times [-\epsilon,0]$, and  $\eta (\partial R_a) \cap R_b$ is a product
neighborhood $\partial R_a\times [0,\epsilon]$.
Assume that $\eta (\partial R_a )\cap {\Cal E'}$ consists of a
collection of
$2r$ parallel vertical annuli.
Let $t$ be a knot  embedded in $R_a \cup R_b$,
such that: $t\cap R_b$ consists of an arc properly embedded in the product
neighborhood $\partial R_a\times [0,\epsilon]$ containing
only a minimum (as defined as in 2.3),  $t\cap R_a$ consists of an arc properly
embedded in the product neighborhood $\partial R_a\times [-\epsilon,0]$,
containing only a maximum,
and $t$ intersects $\Cal E'$ in finitely many points, all lying in
$\eta (\partial R_a)\cap {\Cal E'}$. Suppose these intersection
points lie at different heights (in $\partial R_a\times [-\epsilon,\epsilon]$),
are denoted  by $x_1,x_2,\dots,x_n$, and are ordered according to its height.

\medskip
\noindent {\bf 2.6.1} Let $\Cal E = \Cal E'- int\, \eta (t)$.
We call a solid torus together with a surface $\Cal E$ and an arc
$t$, a piece of type $\tilde \Cal E$. We assume that a piece of type
$\tilde \Cal E$ satisfies the following:

\roster

\item Note that $A_1$ separates a solid torus $N_1 \subset R_a$. A similar
condition as in 2.3.1(4) is satisfied.

\item A similar condition as in 2.3.1(5) is satisfied.

\item Suppose an arc, which does not contain the maximum not the minimum of $t$,
intersects a component of $\Cal E'$,  in two points, say $x_i$, $x_j$, but
does not intersects any other component between these two points. The same
condition as in 2.3.1(3) is satisfied.

\item The subarc of $t$ which contains the minimum has endpoints on
different disks.

\endroster

It is not difficult to construct examples
of pieces satisfying all these conditions.
If the knot $t$ does not intersect the surface $\Cal E'$, then this is a
sphere bounding a 3-ball in a lens space $L(p,q)$ and $t$ is a trivial knot
lying in that ball.

\proclaim{Lemma 2.6.2} The surface $\Cal E$ is incompressible and meridionally
incompressible in $R_a \cup R_b - int\, \eta (t)$.
\endproclaim

\demo{Proof} It is similar to the proof of Lemma 2.4.2.
\enddemo

\remark{Remark} We have a degenerate case of this construction, which happens when
$\Cal E$ is formed by an annulus of meridional slope in $R_a$, and two meridians
disks in $R_b$. In this case $\Cal E$ is a separating sphere in $S^1\times
S^2$ which bounds a 3-ball, and $t$ is a trivial knot in such 3-ball. In this case,
if $t$ intersects $\Cal E$, then clearly the surface will be compressible, i.e.,
a meridian of $R_a$ will make a compression disk.

\medskip
\subhead {2.7 Pieces of type $\tilde \Cal F$} \endsubhead
\medskip

Let $R_a$ and $R_b$ be two solid tori.
Let $A_1,...,A_r$ be a collection of parallel meridian
disk in the solid torus $R_a$ and
let $A_1',...,A_{r}'$ be a collection of parallel meridian disks in a solid torus
$R_b$. Identify $\partial R_a$ with $\partial R_b$, so that  the collection
of curves $\{ \partial A_i\}$ are identified with the curves $\{ \partial A_i'\}$.
Let $\Cal F'$
be the union of the disks $A_i$ and $A_i'$. $\Cal F'$ is a surface  embedded
in $R_a\cup R_b = S^1\times S^2$. Note that each component of $\Cal F'$ is a
non-separating sphere.
Let $\eta(\partial R_a)$ be a
regular neighborhood of $\partial R_a$ in $R_a\cup R_b$ such that
$\eta (\partial R_a) \cap R_a$ is a product neighborhood
$\partial R_a\times [-\epsilon,0]$, and  $\eta (\partial R_a) \cap R_b$ is a product
neighborhood $\partial R_a\times [0,\epsilon]$.
Assume that $\eta (\partial R_a )\cap {\Cal F'}$ consists of a
collection of
$2r$ parallel vertical annuli.
Let $t$ be a knot  embedded in $R_a \cup R_b$,
such that: $t\cap R_b$ consists of an arc properly embedded in the product
neighborhood $\partial R_b\times [0,\epsilon]$ containing
only a minimum (as defined as in 2.3),  $t\cap R_a$ consists of an arc properly
embedded in the product neighborhood $\partial R_a\times [-\epsilon,0]$,
containing only a maximum,
and $t$ intersects $\Cal F'$ in finitely many points, all lying in
$\eta (\partial R_a)\cap {\Cal F'}$. Suppose these intersection
points lie at different heights (in $\partial R_a\times [-\epsilon,\epsilon]$),
are denoted  by $x_1,x_2,\dots,x_n$, and are ordered according to its height.

\medskip

\noindent {\bf 2.7.1} Let $\Cal F = \Cal F'- int\, \eta (t)$.
We call a solid torus together with a surface $\Cal F$ and an arc
$t$, a piece of type $\tilde \Cal F$. We assume that a piece of type
$\tilde \Cal F$ satisfies the following:
\roster

\item Suppose an arc, which does not contain the minimum nor the maximum of $t$,
intersects a component of $\Cal F'$,  in two points, say $x_i$, $x_j$, but
does not intersects any other component between these two points. The same
condition as in 2.3.1(3) is satisfied.

\item The subarc of $t$ which contains the minimum (maximum) has endpoints on
different disks, or in the same disk but intersect it from different sides.

\endroster

It is not difficult to construct examples
of pieces satisfying these conditions.
If the knot $t$ does not intersect the surface $\Cal F'$, then $t$ is a
trivial knot in $S^1 \times S^2$, and if $t$ intersects a component of $\Cal F'$
in one point, then $t$ is of the form $S^1 \times \{* \}$, and any component
of $\Cal F'$ is a sphere of the form $\{*\}\times S^2$.

\proclaim{Lemma 2.7.2} The surface $\Cal F$ is incompressible and meridionally
incompressible in $R_a \cup R_b- int\, \eta (t)$.
\endproclaim

\demo{Proof} It is similar to the proof of Lemma 2.4.2.
\enddemo

\medskip
\subhead {2.8 The meridional surface} \endsubhead
\medskip

\noindent {\bf 2.8.1} Let $T$ be a Heegaard torus in $M$, and let $I=[0,1]$. Consider
$T\times I\subset M$.
$T\times \{ 0\}=T_0$ bounds a solid torus $R_0$, and $T_1=T\times \{1\}$
bounds a
solid torus $R_1$, such that $M= R_0\cup (T\times I)\cup R_1$.
(If $M=S^3$, assume $R_0$ contains the point at infinity). Choose $2n+1$ distinct
points on $I$, $b_0=0,\ a_1, \ b_1,\ a_2,\dots, a_n,\  b_n=1$, so
that $a_i < b_i < a_{i+1}$, for all
$i\leq n-1$. Consider the tori
$T\times \{a_i\}$ and $T\times \{b_i\}$
(so $T\times \{b_0\}= T_0$, $T\times \{b_n\}= T_1$).
If $\beta$ and $\gamma$ are two simple closed curved on $T$, $\Delta(\beta,\gamma)$
denotes, as usual, its minimal geometric intersection number;
if the curves lie on parallel tori, then $\Delta(\beta,\gamma)$ denotes
that number once the curves have been projected into a single torus.

Denote by $\mu_M$ the essential simple closed curve on
$T\times \{b_n\}$ which bounds a meridional disk in $R_1$, and denote by
$\lambda_M$
the essential simple closed curve on $T\times \{b_0\}$ which bounds a disk
in $R_0$. So $\Delta(\mu_M,\lambda_M)=1$ if
$M=S^3$, $\Delta(\mu_M,\lambda_M)=0$ if $M=S^1\times S^2$, and
$\Delta(\mu_M,\lambda_M)=p$ if $M=L(p,q)$.

Let $\gamma_i$ be a simple closed essential curve embedded in the
torus $T\times \{b_i\}$.
Let $\Gamma$ be the collection of all the
curves $\gamma_i$ .

 \medskip

Suppose that  $\Gamma$ satisfies the following:

\roster

\item Either $\gamma_0 =\lambda_M$, or
$\Delta(\gamma_0,\lambda_M) \geq 2$, that is, $\gamma_0$ is not
homotopic to the core of the solid torus $R_0$.

\item $\Delta(\gamma_i,\gamma_{i+1}) \geq 2$, for all $0\leq i\leq n-1$.

\item Either $\gamma_n = \mu_M$, or $\Delta(\gamma_n,\mu_M) \geq 2$,
that is, $\gamma_n$ is not
homotopic to the core of the solid torus $R_1$.

\endroster

Suppose a such $\Gamma = \{ \gamma_i \}$, with $n\geq 2$ is given.
If $\gamma_0$ is a meridian of $R_0$, choose a piece of type $\tilde \Cal C$
in $R_0\cup (T\times [0,a_1])$, denoted $\tilde \Cal C_0$, and which is determined
by a surface $\Cal C_0$ such that
${\Cal C_0}\cap T_0$ is a collection of curves parallel to $\gamma_0$.
If $\gamma_0$ is not a meridian of $R_0$ choose a piece of type
$\tilde \Cal B$ in  $R_0\cup (T\times [0,a_1])$, denoted $\tilde \Cal B_0$,
and determined by a surface $\Cal B_0$
so that ${\Cal B_0}\cap T_0$ is
a collection of curves parallel to $\gamma_0$. For $\gamma_i$, $i\not= 0,n$,
choose a piece of type $\tilde \Cal A$ in $T\times [a_i,a_{i+1}]$, denoted
$\tilde \Cal A_i$, determined by a surface $\Cal A_i$  so that
${\Cal A_i}\cap (T\times \{b_i \})$ is a
collection of curves parallel to $\gamma_i$. If
$\gamma_n$ is a meridian of $R_1$, choose a piece of type $\tilde \Cal C$ in
$R_1\cup (T\times [a_n,1])$, denoted $\tilde \Cal C_n$, determined by a surface
$\Cal C_n$ such that
${\Cal C_n}\cap T_1$ is a collection of curves parallel to $\gamma_n$. If
$\gamma_n$ is not a meridian of $R_1$ choose a piece of type $\tilde \Cal B$
in  $R_1\cup (T\times [a_n,1])$, denoted $\Cal B_n$, determined by a surface
$\Cal B_n$ so that ${\Cal B_n}\cap T_1$ is a collection of curves parallel
to $\gamma_n$. Furthermore suppose that
${\Cal C_0}\cap (T\times\{a_1 \})={\Cal A_{1}}\cap (T\times \{a_1\})$ or
${\Cal B_0}\cap (T\times\{a_1 \})={\Cal A_{1}}\cap (T\times \{a_1\})$,
${\Cal A_i}\cap (T\times\{a_{i+1} \})={\Cal A_{i+1}}\cap (T\times \{a_{i+1}\})$
for $1\leq i \leq n-1$, and
${\Cal C_n}\cap (T\times\{a_n \})={\Cal A_{n}}\cap (T\times \{a_n\})$
or ${\Cal B_n}\cap (T\times\{a_n \})={\Cal A_{n}}\cap (T\times \{a_n\})$, that is,
the boundary curves of the surface in a piece are
identified with the boundary curves of the surface in an adjacent piece. Suppose also
that the endpoints of the arcs coincide so that the union of all the arcs is a
knot $K$ in $M$. Note
that the union of all the surfaces becomes a surface $S$ properly embedded in
$M-int\, \eta (K)$, whose boundary consists of meridians of the knot $K$, so
$S$ is a meridional surface. It is not difficult to construct examples of surfaces
satisfying  those conditions.

Suppose now $\Gamma =\{\gamma_0\}$. In
this case the surface is given by a piece of type $\tilde \Cal D$, $\tilde
\Cal E$ or $\tilde \Cal F$. More precisely, if $\gamma_0$ is not a meridian
of $R_0$ and not a meridian of $R_1$, then take a piece of type $\tilde \Cal D$,
so that the intersection of the surface $\Cal D$ with the torus $T\times \{b_0\}$
consists of curves parallel to $\gamma_0$.
Similarly, if $\gamma_0$ is a meridian of  only one of $R_0$ or $R_1$, take a piece
of type $\tilde \Cal E$, and if $\gamma_0$ is a meridian of both solid tori,
then take a piece of type $\tilde \Cal F$, or a degenerate piece of type
$\tilde \Cal E$.

\proclaim {Theorem 2.8.2} Let $K$ be a knot and $S$ a meridional surface as
constructed in 2.8.1.
Then $K$ is a $(1,1)$-knot and $S$ is an essential meridional surface.
\endproclaim

\demo {Proof} As in 2.8.1, let  $T$ be a Heegaard torus in $M$, and let
$I=[0,1]$. Consider $T\times I\subset M$. $T\times \{ 0\}=T_0$ bounds a
solid torus $R_0$, and $T_1=T\times \{1\}$ bounds a solid torus $R_1$, such
that $M= R_0\cup (T\times I)\cup R_1$. As in 2.8.1, there are $2n+1$ distinct
points on $I$, $b_0=0,\ a_1, \ b_1,\ a_2,\dots, a_n,\  b_n=1$, so that $a_i
< b_i < a_{i+1}$, for all
$i\leq n-1$.

If $S$ is made of one piece, i.e., it is of type $\tilde \Cal D$,
$\tilde \Cal E$ or $\tilde \Cal F$, then it follows from Lemmas
2.5.2, 2.6.2 and 2.7.2 that the surface is essential.
Suppose then that the $S$ is made of several pieces.

In each level $a_i$, consider the nested disks $D_i$, and between each
two consecutive levels $a_i, a_{i+1}$, take the annuli $B,B',C,C'$ as
defined 2.2.1, 2.3.1 or 2.4.1.
Suppose $S$ is compressible in $S^3-K$ or meridionally compressible, and
let $D$ be a compression disk which intersects $K$ in at most one point.
Assume that $D$ intersects transversely the annuli $B,B',C,C'$ and the
disks $D_i$. Then the intersection consists of a finite collection of simple
closed curves and arcs. Simple closed curves of intersection can be removed
as in Lemmas 2.2.2, 2.3.2, 2.4.2. Suppose then that the intersections between
$D$ and the annuli
and disks consists only of arcs. Let $\alpha$ be a such outermost arc. This
can be eliminated as in those lemmas, except in two cases.

First, the arc $\alpha$ lies in a disk $D_1$, and the disk $D' \subset D$ cutoff by $\alpha$
lies in a region between $D_1$ and an annulus $B'$ or $C'$. It follows that one of
the conditions 2.2.1(1), 2.3.1(1) or 2.4.1(1) is not satisfied

Second, the disk $D' \subset D$ cutoff by $\alpha$ lies in a region
between two annuli of type $B$ or  $C$, i.e., a region where two different
pieces of type $\tilde \Cal A$ are glued, or where $\tilde \Cal C_0$ or
$\tilde \Cal B_0$ are glued to $\tilde \Cal A_1$. The boundary of this disk will consist
of an arc on lying on $S$ and the arc $\alpha$, which lies on $B$.
Now condition 2.8.1(2) and Lemma 3.1 of [E1]
show that this is not possible.
Then the disk $D$ is disjoint from the collection of annuli
and disks.  Again, by the lemmas, the only possibility left, is that the disk
lies in a region  between annuli of type $B$ and $C$, and then it is disjoint
from $K$. Again by Condition 2.8.1(2) and Lemma 3.1 of [E1] this is not
possible.

Note that by construction the knot $K$ lies in a product
$T\times [-\epsilon,1+\epsilon]$, so that it has a maximum at
$T\times \{-\epsilon \}$, a minimum at $T\times \{1+\epsilon \}$, and intersect any
other horizontal torus in 2 points. From this follows that $K$ is a $(1,1)$-knot.
\qed\enddemo

Note that if a surface $S$ does not satisfy one of the conditions 2.8.1(1)-(3), then it
will be compressible. For example, suppose that condition 2.8.1(1) is not satisfied.
Then the piece $\tilde \Cal B_0$ is made with annuli which go
once longitudinally around $R_0$. If $t$ intersect $\Cal B_0$,
it was observed in 2.3.1 that the surface will be compressible;
if the knot is disjoint from $\Cal B_0$, then this part of the surface $S$ is
isotopic into the level torus $T\times \{a_1\}$, and then there is a compression disk
lying in $T\times [a_1,b_1]$.

It follows from this construction that for any $M$, and any given
integers $g\geq 1$ and $h\geq 0$, there exist $(1,1)$-knots which
admit a meridional essential surface of genus $g$ and $2h$ boundary
components. In fact, there are knots which contain two or more
meridional surfaces which are not parallel in the knot complement.
Note that if in the construction the knot is disjoint from the
surface, then what we get is one of the knots and surfaces
constructed in [E1], possibly with several parallel components. As
we said before, $(1,1)$-knots in $S^3$ do not admit any meridional
essential  surface of genus $0$, but when $M=L(p,q)$, it follows
that for any given integer $h\geq 1$, there exist $(1,1)$-knots
which admit a meridional essential surface of genus $0$ and $2h$
boundary  components. These are given by pieces of type $\tilde \Cal
E$. In particular, there exists composite $(1,1)$-knots in lens
spaces, i.e., knots which admit an essential meridional planar
surface with two boundary components. It follows from 3.6.1(1) that
these knots are obtained as a connected sum of the core of a
Heegaard torus in a lens space with a two-bridge knot in the
3-sphere.


As a special case consider $(1,1)$-knots which admit a meridional surface of genus 1.
Note that a piece of type $\tilde \Cal A$ or of type $\tilde \Cal B$ consists
of punctured tori, so a surface $S$ of genus 1 constructed as in 2.8.1 cannot contain
two pieces of type $\tilde \Cal A$, or two pieces of type $\tilde \Cal B$,
or a piece of type $\tilde \Cal A$ and a piece of type $\tilde \Cal B$.
So a surface of genus 1 must consists either of a piece of type $\tilde \Cal A$ and two
pieces of type $\tilde \Cal C$, or of a piece of type $\Cal B$ and a piece of type
$\tilde \Cal C$, or it is just a piece of type $\tilde \Cal D$.
Suppose now we have a meridional surface of genus 1 with 2 punctures in the exterior
of a $(1,1)$-knot. Note that in a piece of type $\tilde \Cal C$, each component always
intersects the knot in at least two meridians. So if we have a meridional surface of genus 1
with two punctures, then it must consist either of a piece of type $\tilde \Cal B$ and a
piece of type $\tilde \Cal C$, or it is a piece of type $\tilde \Cal D$.
In the first case, the piece of type $\tilde \Cal B$ consist of just a once punctured torus
without intersections with the knot, and the piece of type $\tilde \Cal C$ consist of a
disk intersecting the knot in two meridians, so the surface looks like a torus with a bubble.
In the second case, the piece of type $\tilde \Cal D$ consists of a torus which either
is the boundary of a regular neighborhood of a torus knot in $M$, or it is isotopic (in $M$)
to the Heegaard torus $T$, as in Figure 5.

It follows from [CGLS, 2.0.3] that any of the knots constructed here will have a
closed incompressible surface in its complement, which will be in general meridionally
compressible. However, we do not intend here to determine which
is that surface.

\medskip
\head {3. Characterization of meridional surfaces} \endhead
\medskip

\proclaim {Theorem 3.1} Let $K$ be a $(1,1)$-knot in a closed 3-manifold $M$,
 and let $S$ be an essential
meridional surface for $K$. Then $K$ and $S$ come
from the construction of 2.8, that is, $K$ and $S$ can be isotoped so that
they look as one of the knots and surfaces constructed in 2.8.
\endproclaim

\demo{Proof} Let $T$ be a Heegaard torus in $M$, and let $I=[0,1]$. Consider
$T\times I\subset M$.
$T\times \{ 0\}=T_0$ bounds a solid torus $R_0$, and $T_1=T\times \{1\}$
bounds a solid torus $R_1$, such that $M= R_0\cup (T\times I)\cup R_1$. Let $k$ be a
$(1,1)$-knot,  and assume that $k$ lies in $T\times I$, such that
$k\cap T\times \{ 0\} = k_0$ is an  arc,
$k\cap T\times \{ 1\} = k_1$ is an arc, and $k\cap T\times (0,1)$ consists
of two straight arcs.

Suppose there is a meridional surface $S$ in $M-int\, \eta (k)$, which is incompressible
and meridionally  incompressible. Consider $S$ as a surface embedded in $M$ which
intersects $k$ in a finite number of points. Assume that
$S$ intersects transversely $T_0$ and
$T_1$. Let $S_0=S\cap R_0$, $S_1=S\cap R_1$, and
$\tilde S=S\cap (T\times I)$. We can assume that all the intersection points
between $S$ and $k$ lie on $\tilde S$.
Let $\pi:T\times I \rightarrow I$ be the height function, where we
choose $0$ to be
the highest point, and $1$ the lowest. We may assume that the height
function on
$\tilde S$ is a Morse function. So there is a finite set of different points
$X=\{x_1,x_2,\dots,x_m\}$ in $I$, so that
$\tilde S$ is tangent to $T\times \{ x_i \}$ at exactly one point, and this
singularity can be a local maximum, a local minimum, or a simple saddle.
Assume also that no point of intersection between $\tilde S$ and $k$ lie
on one of the levels  $T\times \{ x_i \}$.
For any $y\notin X$,
$T\times \{ y \}$ intersects $\tilde S$ transversely, so for any such $y$,
$\tilde S \cap T\times \{ y \}$ consists of a finite collection of simple
closed
curves called level curves, and at a saddle point $x_i$, either one level
curve of
$\tilde S$ splits into two level curves, or two level curves are fused
into one curve.

Define the complexity of $S$ by the pair
$c(S)=(\vert S_0\vert + \vert S_1 \vert + \vert \tilde S \vert , \vert X \vert)$
(where $\vert Y \vert$ denotes the number of points if $Y$ is a finite set,
or the number of connected components if it is a surface, and give to such
pairs the lexicographical order). Assume that $S$ has been isotoped so that
$c(S)$ is minimal.


\proclaim {Claim 3.2} The surfaces $S_0$, $S_1$ and $\tilde S$ are
incompressible and meridionally incompressible in $R_0$, $R_1$, and
$T\times I - int\, \eta (k)$ respectively.
\endproclaim

\demo{Proof}  Suppose one of the surfaces is compressible or meridionally
compressible, say $\tilde S$, and let
$D$ be a compression disk, which it is disjoint from $k$, or intersects it
in one point. Then $\partial D$ is
essential in $\tilde S$ but inessential in $S$. By cutting $S$ along $D$ we
get a surface $S'$ and a sphere $E$. Note that $S$ and $S'$ are isotopic in
$M-k$ (see the remark below the proof of this claim).
For $S'$ we can similarly define the surfaces
$S_0'$, $S_1'$ and
$\tilde S'$. Note that
$\vert S_i \vert = \vert S_i' \vert + \vert E \cap R_i
\vert$, $i=1,\ 0$, then either $\vert S_0' \vert < \vert S_0 \vert$ or
$\vert S_1' \vert < \vert S_1 \vert$, for $E$ intersects at least one of
$R_0$,
$R_1$. Also $\vert \tilde S' \vert \leq \vert \tilde S \vert$, so
$c(S') < c(S)$, but this contradicts the choice of $S$.
\qed
\enddemo

\remark {Remark} The surfaces $S$ and $S'$ could be non-isotopic if one
of the following cases occur: a) There is a non-separating sphere disjoint
from the knot; or b) there is a separating sphere intersecting the knot
in two points, which bounds a 3-ball containing a non-prime knotted arc
of the knot. In the first case it is not difficult to see that the knot
has to be trivial, and it follows from the proof of
Theorem 3.1 that the second case is impossible, for any such sphere
determines a piece of type $\tilde \Cal E$ which has a prime summand,
in fact a 2-bridge knot summand.
\endremark

This implies that $S_0$ is a collection of trivial disks, meridian
disks and essential annuli in $R_0$. If a component of $S_0$ is a trivial
disk $E$,
then $\partial E$ bounds a disk on $T_0$ which contains $k_0$, for
otherwise
$\vert S_0 \vert$ could be reduced. If a component of $S_0$ is an essential
annulus
$A$, then $A$ is parallel to an annulus $A^\prime \subset T_0$, and
$A^\prime$ must contain $k_0$, for otherwise $\vert S_0 \vert$ could be
reduced.
This implies that the slope of $\partial A$ cannot consists of one
meridian and
several longitudes, for in this case
$A$ would also be parallel to $T_0-A^\prime$, and then $\vert S_0\vert$
could be reduced. This also implies that $S_0$ cannot contain both
essential annuli and meridian disks. A similar thing can be said for $S_1$.

\proclaim {Claim 3.3} $\tilde S$ does not have any local maximum or minimum.
\endproclaim

\demo{Proof} Suppose $\tilde S$ does have a maximum. It will be shown that
$\tilde S$ has a component which is either a disk, an annulus, a once
punctured
annulus, or a once punctured torus which is parallel to a subsurface in
$T_1$.
Choose the maximum at lowest level, say at level
$x_i$, so that there are no other maxima between $x_i$ and 1. So
$\tilde S\cap (T\times \{ x_i \})$ consists of a point and a collection of
simple
closed curves. Just below the level $\{x_i\}$, the surface $\tilde S$
intersects
the level tori in simple closed curves, so below the maximum, a disk $E_1$
is
being formed; if $x_i$ is the last singular point then the disk $E_1$ will
be
parallel to a disk on $T_1$.

Look at the next singular level
$x_{i+1}$, $x_i<x_{i+1}$. Note that the disk $E_1$ may intersect $k$. If
the singular point at $x_{i+1}$ is a local minimum, or a saddle whose level
curves are  disjoint from the boundary curve of $E_1$, then we can
interchange the singularities.  If it is a saddle point involving the curve
$\partial E_1$ and another curve, then both singular points cancel each
other, by pushing down the maximum; this can also be done if one of the arcs
of $k$ intersects $E_1$, for a given arc cannot intersect $E_1$ more than
once (for otherwise a subarc of $k$ will be isotopic to an arc on $S$,
and then $S$  would be compressible). If the singular point at
$x_{i+1}$ is a saddle formed by a selfintersection of the curve $\partial
E_1$, then the disk $E_1$ below the maximum transforms into an annulus
$E_2$. If this singularity occurs inside a 3-ball bounded by $E_1$ and
a disk of a level torus, then it is not difficult to see that $E_2$ would
be compressible or meridionally compressible in $T\times [y_1,y_2]$, for
certain levels $y_1$, $y_2$, and then either $\tilde S$ would be compressible or
meridionally compressible or could be isotoped to reduce the number of
singular points. So the singularity occurs outside such 3-ball, and then the
annulus $E_2$  must be parallel to an annulus in some $T\times \{y\}$.
For any $y$ just below $x_{i+1}$, $\partial E_2$ consists of two parallel
curves $c_1$ and $c_2$. If $c_1$ is a trivial curve on
$T\times \{ y\}$ which bounds a disk disjoint from $k$, then by cutting
$\tilde S$ with such disk, we get a surface $S'$ isotopic to $S$, but with
$c(S')< c(S)$, for
$S'$ has less singular points than $S$. If $c_1$ bounds a disk which
intersects
$k$ once, then $S$ would be meridionally compressible. So there remain two
possibilities, either $c_1$ and $c_2$ are essential curves on
$T\times \{ y\}$ or $c_1$ bounds a disk which intersects $k$ twice. Note that
if $k$ intersects $E_1$ then necessarily $c_1$ and $c_2$ are essential
curves on $T\times \{ y\}$.  Look at the next singular level $x_{i+2}$. If
it is a local minimum or a saddle whose level curves are disjoint from the
annulus $E_2$, then again we can interchange the singularities (note that
these curves cannot be inside the solid torus bounded by
$E_2$ and an annulus in  $T\times \{ x_{i+2} \}$, for there would be a
maximum in there). If it is a saddle joining the annulus $E_2$ with another
curve, then by pushing $E_2$ down the singular point $x_{i+2}$ is eliminated.
Again note that this can also be done if one of the arcs of $k$
intersects $E_2$, for a given arc cannot intersect $E_2$ more than once, for
otherwise $S$ would be compressible. So this singular point has to be of the
annulus with itself, and then a new surface $E_3$ is formed. The surface
$E_3$ has to be a once punctured torus or a once punctured annulus. Suppose
first that the surface $E_3$ is not parallel to some $T\times\{y^\prime\}$,
i.e., the singular point $x_{i+2}$ occurs inside the solid torus bounded by
$E_2$ and an annulus in  $T\times \{ x_{i+2} \}$. It is not difficult to see
that in this case the surface $E_3$ is compressible or meridionally
compressible; this is because the arcs of $k$ are straight, and if they look
complicated inside $E_2$, we can find a level preserving isotopy which make
them look ``straight''. Then the surface $E_3$ is parallel to some surface in
$T\times\{y^\prime\}$. Note that if the knot intersects $E_2$ then this
surface will be compressible or meridionally compressible. If not, then
there may be one more singular point of $E_3$ with itself, but then the
surface
$E_4$ which will be formed would be compressible, for a component of its
boundary would bound a disk in some $T\times \{y''\}$ disjoint from $k$.

We conclude that either $\tilde S$ is compressible, meridionally compressible,
the number of singular points is not minimal, or a
component of $\tilde S$, say
$S^\prime$, it is parallel to a surface in $T_1$. In the latter case, if the
surface $S^\prime$ is disjoint from $k$, then it can be pushed to $R_1$,
reducing $c(S)$. If it not disjoint from $k$, then by pushing $S'$ into
$T_1$, we see that the arc $k_1$ would be
parallel to $S$, implying that $S$ is compressible.
\qed
\enddemo

At a nonsingular level $y$, $\tilde S\cap (T\times \{ y \})$ consists of
simple closed
curves. If such a curve $\gamma$ is trivial in $T\times \{ y \}$, then it
bounds a
disk in that torus.  If such a disk is disjoint from $k$, then by
the incompressibility of $\tilde S$, $\gamma$ bounds a disk in $\tilde S$,
which give
rise to a single maximum or minimum, contradicting Claim 3.3. Such a disk
cannot
intersect $k$ once, for
$S$ is meridionally incompressible, so $\gamma$ bounds a disk which
intersects $k$
twice.

\proclaim {Claim 3.4} Only the following types of saddle points are possible.
\roster
\item A saddle changing a trivial simple closed curve into two essential
simple
closed curves.
\item A saddle changing two parallel essential curves into a trivial curve.
\endroster
\endproclaim

\demo{Proof} At a saddle, either one level curve of $\tilde S$ splits
into two level curves, or two level curves are joined into one level curve.
If a level curve is trivial in the corresponding level torus and at a saddle
the curve joins with itself, then the result must be two essential simple
closed curves, for if the curves obtained are trivial, then $\tilde S$ would
be
compressible or meridionally compressible. If a curve is nontrivial and at
the
saddle joins with itself, then the result is a curve with the same slope
as the
original and a trivial curve, for the saddle must join points on the same
side of
the curve, by orientability.  If two trivial level curves are joined into
one, then
because the curves must be concentric, the surface $S$ would be
compressible or meridionally compressible.
So only the following types of saddle points are possible:

\roster
\item A saddle changing a trivial simple closed curve into two essential
simple
closed curves.
\item A saddle changing two parallel essential curves into a trivial curve.
\item A saddle changing an essential curve $\gamma$ into a curve with the
same slope as $\gamma$, and a trivial curve.
\item A saddle changing an essential curve $\gamma$ and a trivial curve
into an
essential curve with the same slope as $\gamma$.
\endroster

We want to show that only saddles of types 1 and 2 are possible. Look at
two
consecutive saddle points. We note that in many cases $\tilde S$ can be
isotoped
so that two consecutive singular
points can be put in the same level, and in that case a compression disk for
$\tilde S$ can be found. Namely,

\roster
\item A type 3 is followed by a type 1, see Figure 6(a).
\item A type 3 followed by a type 2, see Figure 6(b).
\item A type 1 is followed by a type 4, see Figure 6(c).
\item A type 2 followed by a type 4, see Figure 6(d).
\item A type 3 followed by a type 4, see Figure 6(e) and 6(f).
\endroster

\midinsert
\figure{2.5}{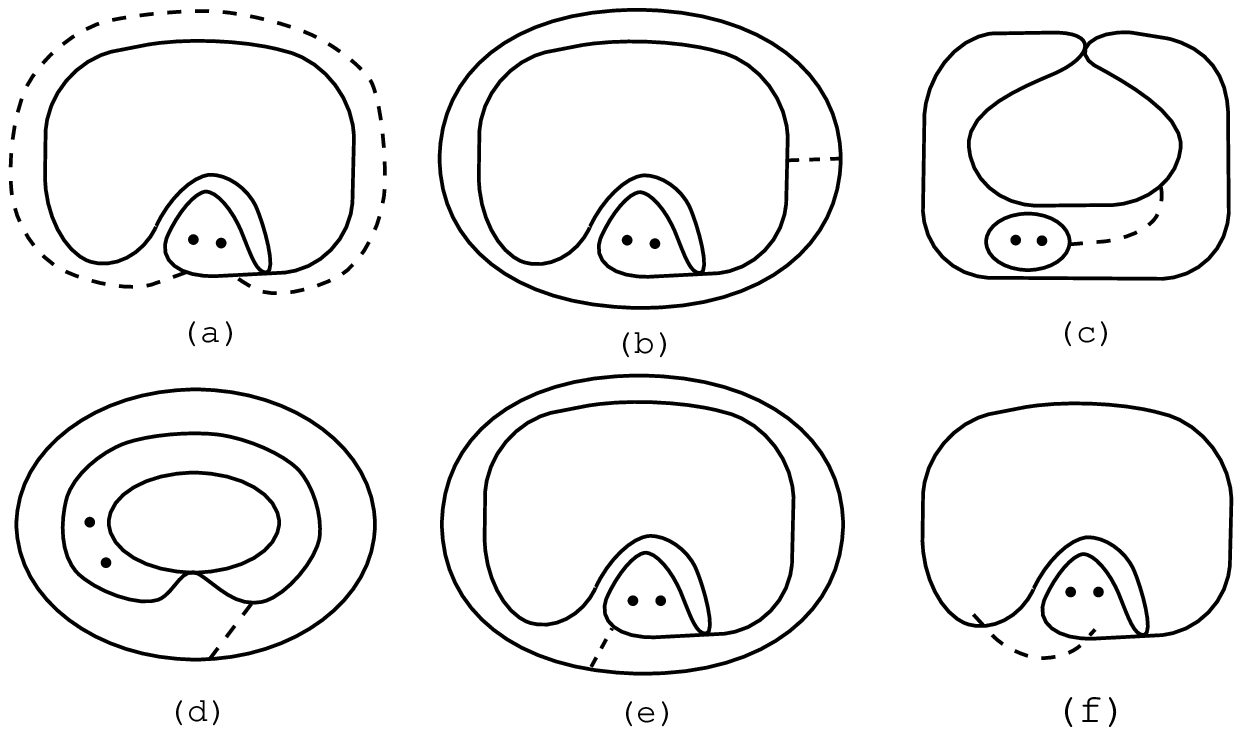}{6}
\endinsert

The dots in Figure 6 indicate the points of intersection between $k$ and
some
$T\times \{ y \}$, and the dotted line indicates the path followed by the
second
saddle point. Note that in all the other cases it may not be possible to
put the
two saddles points in the same level. Note also that the knot does not
intersect $\tilde S$ between two of these singular points in any of the six
cases, because there is a trivial curve in each level bounding a disk which
intersects $k$ twice.

Suppose there are singularities of type 3, and take the lowest one, say at
level $x_m$. If there is a singularity at level $x_{m+1}$, then it is of type
1, 2, or 4,
which contradicts the previous observation. So $x_m$ is in fact the last
singularity. Again note that there is no intersections of $k$ with the
surface $\tilde S$ after the singular point $x_m$. Then $S_1$ consists of at
least a disk, and at least one annulus or a meridian disk. Isotope $S$ so
that the saddle at $x_m$ is pushed into $R_1$, this can be done without
moving $k_1$. If
$S_1$ contained a meridian disk and a trivial disk, which joined in $x_m$,
then now there is just one meridian disk, reducing $\vert S_1\vert$ and
$\vert X\vert$ by 1, then reducing $c(S)$. If $S_1$ contained
an annulus and a disk, which joined in the singularity, then now there is
just one annulus, again reducing $c(S)$.

If there is a saddle of type 4, then take the highest one. If it is not at
level $x_1$, then it is preceded by a singularity of type 1, 2, or 3, which is
a contradiction. So it is a level $x_1$. Push the singularity to
$R_0$. Doing an argument as in the previous case we get a contradiction.
\qed
\enddemo

\proclaim {Claim 3.5} $S_0$ ($S_1$) consists only of annuli or only of
meridian disks.
\endproclaim
\demo{Proof} Suppose $S_0$ contains a trivial disk. Look at the first
singularity.
Suppose it is of type 2. If there is an annulus in $S_0$, then its boundary
components are joined in the saddle, but because there is the trivial disk,
it is not difficult to see this implies that $S$ is compressible.
If there is
a pair of meridian disks in $S_0$, and their boundaries are joined in the
saddle, then we get a disk which can be pushed to $S_0$, reducing $\vert S_0
\vert$ by 1. So the first singularity has to be of type 1, and a trivial
disk touches with itself. Push the singularity to $R_0$, so instead of the
disk we have now an annulus, this leaves $\vert S_0\vert$ unchanged, but
reduces $\vert X\vert$ by one.
\qed
\enddemo

\proclaim {Claim 3.6} If $\vert X \vert =0$, then $S$ is a piece of type
$\tilde \Cal D$, $\tilde \Cal E$ or $\tilde \Cal F$.
\endproclaim

\demo{Proof} As $\tilde S$ has no singular points, then it has to be a
collection of annuli. Then $S_0$ and $S_1$ determine the same slope in
$T_0$ and $T_1$ respectively. If $S_0$ and $S_1$ consist of annuli then we
have a piece of type $\tilde \Cal D$, if one consists of annuli and the
other of meridian disks then we have a piece of type $\tilde \Cal E$, and if
both consist of meridian disks then we have a piece of type $\tilde \Cal F$.
Note that the conditions 2.5.1, 2.6.1 or 2.7.1 for the arcs must be satisfied for
otherwise the surface $S$ would be compressible or meridionally compressible.
\qed
\enddemo

\proclaim {Claim 3.7} Suppose that $\vert X \vert > 0$. Then $S$ is a meridional
surface as in 2.8, and it is made of a union of pieces of type $\tilde \Cal
A$, $\tilde \Cal B$ and $\tilde \Cal C$.
\endproclaim

\demo{Proof} The surface $S_0$ consists only of annuli or only of meridian disks, so
$\partial S_0$ consists of $p$ essential curves. This implies that the first
saddle is of type 2, changing two essential curves into a trivial curve. The
knot $k$ may intersect the surface before the first singularity, but it
cannot intersect it right after that singularity, for a trivial curve is
formed. If $p>2$, then the next singularity is again of type 2, for if it is
type 1, the new essential curves that are formed will have the same slope as
the original curves, and then it is not difficult to see that there is a
compression disk for $\tilde S$. So there will be singularities of type 2,
until no essential curve is left, so there are in total $p/2$ of these
singularities. In particular this shows that $p$ is even; this was clear if
$S_0$ consists of annuli, but it is not clear if it consists of meridian
disks. This shows that $S_0$ and a part of $\tilde S$ form a piece of type
$\tilde \Cal B$ or type $\tilde \Cal C$, depending if $S_0$ consists of
annuli or meridian disks.  After the $p/2$ singularities, the intersection
of $S_0$ with a level torus consists of $p/2$ trivial curves, which are nested,
and bound a disk which intersect $k$ in two points. The next singularity has
to be of type 1, changing a trivial curve into 2 essential curves. If the next
singularity after that one is of type 2, again we have a compression disk, so the
singularities are all of type 1, until there is no more trivial curves. Then
there are $p/2$ of such singularities. The knot may intersect the surface
only after the $p/2$ singularities of type $1$, i.e., after there are no
trivial curves in a level, for otherwise $S$ would be compressible or meridionally
compressible. If there are no more singularities, then we get at level $1$, and then
$S_1$ consists of annuli or meridian disks, and then we have another piece of type
$\tilde \Cal B$ or $\tilde \Cal C$. If there are more singularities, the next $p/2$
singularities are of type 2, and then the next ones are of type 1. Note again that the
knot cannot intersect the surface in the levels in which there are trivial curves of
intersection. This will form a piece of type $\tilde \Cal A$.  Continuing in this way,
it follows that $S$ is made of pieces of type $\tilde \Cal A$, $\tilde \Cal B$ or
$\tilde \Cal C$. The conditions on the curves and on the arcs must be satisfied, for
otherwise the surface will be compressible or meridionally compressible.
\qed
\enddemo

This completes the proof of Theorem 3.1.
\qed
\enddemo

\vskip20pt

\noindent {\bf Acknowledgment.}
Part of this work was done while we were visiting the
Department of Mathematics of The University of Texas at Austin. We are
grateful for their hospitality. The first author was supported by CONACYT grant
36507-E and by PAPIIT-UNAM grant IN115105.

\head {References} \endhead

\bigskip

\noindent   {Instituto de Matem\'aticas, UNAM, Ciudad Universitaria,
04510 M\'exico D.F.,\break  MEXICO}

\noindent E-mail: mario\@ matem.unam.mx
\medskip
\noindent { Centro de Investigaci\'on en Matem\'aticas, Apdo. Postal 402,
36000 Guanajuato, Gto., MEXICO}

\noindent E-mail: kikis\@ cimat.mx

\end
\widestnumber\key{CGLS}
\Refs\nofrills{}


\ref\key CM
\by A. Cattabriga and M. Mulazzani
\paper $(1,1)$-knots via the mapping class group of the twice punctured torus
\jour Adv. Geom.  \vol 4  \yr 2004 \pages 263-277
\endref

\ref\key CK
\by Doo Ho Choi and Ki Hyoung Ko
\paper Parameterizations of 1-bridge torus knots,
\jour  J. Knot Theory Ramifications
\vol 12  \year 2003  \pages 463-491
\endref

\ref\key CGLS
\by M. Culler, C. McA. Gordon, J. Luecke and P.B. Shalen
\paper Dehn surgery on knots
\jour Annals of Mathematics
\vol 125 \year 1987 \pages 237-300
\endref

\ref\key E1
\by M. Eudave-Mu\~noz
\paper Incompressible surfaces in tunnel number one knot complements
\jour Topology Appl.
\vol 98 \yr 1999 \pages 167-189
\endref

\ref\key E2
\bysame
\paper Meridional essential surfaces for tunnel number one knots
\jour Bol. Soc. Mat. Mex. (3)
\vol 6 \yr 2000 \pages 263-277
\endref

\ref\key E3
\bysame
\paper Incompressible surfaces and $(1,1)$-knots
\jour to appear in J. Knot Theory Ramifications
\vol \yr  \pages
\endref


\ref\key F
\by E. Finkelstein
\paper Closed incompressible surfaces in closed braid complements
\jour J. Knot Theory Ramifications  \vol 7  \yr 1998 \pages 335-379
\endref

\ref\key GMM
\by  H. Goda, H. Matsuda and T. Morifuji
\paper Knot Floer homology of (1,1)-knots
\jour  Geom. Dedicata
\vol 112 \yr 2005 \pages 197-214
\endref

\ref\key GL
\by C. McA. Gordon and R.A. Litherland
\paper Incompressible surfaces in branched coverings
\inbook The Smith conjecture (New York, 1979), Pure Appl. Math., 112
\pages 139-152
\endref

\ref\key GR
\by C.McA. Gordon and A. Reid
\paper Tangle decompositions of tunnel number one knots and links
\jour J. Knot Theory Ramifications
\vol 4 \yr 1995 \pages 389-409
\endref

\ref\key HT
\by A, Hatcher and W. Thurston
\paper   Incompressible surfaces in $2$-bridge knot complements
\jour Invent. Math.
\vol 79 \yr 1985 \pages 225-246
\endref

\ref\key LP
\by M.T. Lozano and J. Przytycki
\paper Incompressible surfaces in the exterior of a closed $3$-braid. I.
Surfaces with horizontal boundary components.
\jour Math. Proc. Cambridge Philos. Soc.  \vol 98  \yr 1985  \pages 275-299
\endref

\ref\key M
\by W. Menasco
\paper Closed incompressible surfaces in alternating knot and link
complements
\jour Topology
\vol 23 \yr 1984 \pages 37-44
\endref


\ref\key MS
\by K. Morimoto and M. Sakuma
\paper On unknotting tunnels for knots
\jour Math. Ann. \vol 289 \yr 1991 \pages 143-167
\endref


\ref\key O
\by U. Oertel
\paper  Closed incompressible surfaces in complements of star links
\jour Pacific J. Math.  \vol 111  \yr 1984 \pages 209-230
\endref


\endRefs
